\documentclass[11pt, a4paper]{amsart}
\usepackage{amssymb}
\usepackage{mathrsfs}
\usepackage[english]{babel}
\swapnumbers

\begin{document}
\newtheoremstyle{mytheorem}
  {}
  {}
  {\itshape}
  {}
  {\scshape}
  {.}
  { }
  {}

\newtheoremstyle{mydefinition}
  {}
  {}
  {\upshape}
  {}
  {\scshape}
  {.}
  { }
  {}

\newtheoremstyle{theoremappendix}
  {}
  {}
  {\itshape}
  {}
  {\scshape}
  {.}
  { }
  {}

\newtheoremstyle{definitionappendix}
  {}
  {}
  {\upshape}
  {}
  {\scshape}
  {.}
  { }
  {}

%
%
\theoremstyle{theoremappendix}
\newtheorem{thmapp}{Theorem}
\newtheorem{lemmaapp}[thmapp]{Lemma}
\theoremstyle{definitionappendix}
\newtheorem{remapp}[thmapp]{Remark}
\newtheorem{defiapp}[thmapp]{Definition}
\theoremstyle{mytheorem}
\newtheorem{lemma}{Lemma}
\newtheorem*{lemma*}{Lemma}
\newtheorem{prop}[lemma]{Proposition}
\newtheorem{thm}[lemma]{Theorem}
\newtheorem{addendum}[lemma]{Addendum}
\newtheorem*{thm*}{Theorem}
\newtheorem{cor}[lemma]{Corollary}
\newtheorem*{cor*}{Corollary}
\theoremstyle{mydefinition}
\newtheorem{rem}[lemma]{Remark}
\newtheorem*{rem*}{Remark}
\newtheorem{conj}[lemma]{Conjecture}
\newtheorem{notation}[lemma]{Notation}
\newtheorem*{notation*}{Notation}
\newtheorem*{warning*}{Warning}
\newtheorem{rems}[lemma]{Remarks}
\newtheorem{defi}[lemma]{Definition}
\newtheorem*{defi*}{Definition}
\newtheorem{defis}[lemma]{Definitions}
\newtheorem{exo}[lemma]{Example}
\newtheorem{exos}[lemma]{Examples}
\def\isom{\mathrm{Is}}
\def\stab{\mathrm{Stab}}
\def\homeo{\mathrm{Homeo}}
\def\aut{\mathrm{Aut}}
\def\supp{\mathrm{supp}}
\def\esssup{\mathop{\mathrm{ess\,sup}}}
\def\clap#1{\hbox to 0pt{\hss#1\hss}}
\def\mathclap{\mathpalette\mathclapinternal}
\def\mathclapinternal#1#2{%
	\clap{$\mathsurround=0pt#1{#2}$}}
\def\rprod{{^\mathrm{R}}\!\!\mathop{\prod}}
\def\cvx#1{{{C}(#1)}}
\def\A{\mathscr{A}}
\def\B{\mathscr{B}}
\def\C{\mathscr{C}}
\def\D{\mathscr{D}}
\def\E{\mathscr{E}}
\def\F{\mathscr{F}}
\def\G{\mathscr{G}}
\def\M{\mathscr{M}}
\def\O{\mathscr{O}}
\def\T{\mathscr{T}}
\def\V{\mathscr{V}}
\def\W{\mathscr{W}}
\def\geod{\mathfrak{G}}
\def\rays{\mathfrak{R}}
\def\proba{\mathcal P}
\def\ball{\mathsf{B}}
\def\cball{\overline\ball}
\def\one{\mathbf{1\kern-1.6mm 1}}
\def\cat#1{\ensuremath{\mathrm{CAT}(#1)}}
\def\id{{\it I\! d}}
\def\bbb#1{\overline{\!#1}}
\def\tc{\ensuremath{\T_\mathrm{c}}}
\def\tw{\ensuremath{\T_\mathrm{w}}}
\def\weak{weak-* }
\def\se{\subseteq}
\def\p{\partial}
\def\d{\,\mathrm{d}}
\def\ro{\varrho}
\def\fhi{\varphi}
\def\teta{\vartheta}
\def\epsi{\varepsilon}
\def\wt{\widetilde}
\def\defq{\stackrel{\text{\rm\tiny def}}{=}}
\def\ti{-\allowhyphens}
\def\lra{\longrightarrow}
\def\No{N\raise4pt\hbox{\tiny o}\kern+.2em}
\def\no{n\raise4pt\hbox{\tiny o}\kern+.2em}
\def\bsl{\backslash}
\def\beq{\begin{equation}}
\def\eeq{\end{equation}}
%


\title[Superrigidity and splitting]{Superrigidity for irreducible lattices\\ and geometric splitting}
\author{Nicolas Monod}
\address{University of Chicago}
\begin{abstract}
We prove general superrigidity results for actions of irreducible lattices on \cat0 spaces; first, in terms of the ideal boundary, and then for the intrinsic geometry (including for infinite-dimensional spaces). In particular, one obtains a new and self-contained proof of Margulis' superrigidity theorem for uniform irreducible lattices in non-simple groups. The proofs rely on simple geometric arguments, including a splitting theorem which can be viewed as an infinite-dimensional (and singular) generalization of the Lawson-Yau/Gromoll-Wolf theorem. Appendix~A gives a very elementary proof of commensurator superrigidity; Appendix~B proves that all our results also hold for certain non-uniform lattices.
\end{abstract}
\maketitle

\section{Introduction}
\label{sec_intro}%
\subsection{Superrigidity} In the early seventies, Margulis proved his celebrated superrigidity theorem for irreducible lattices in semi-simple Lie and algebraic groups of higher rank. One of the motivations for this result is that it implies \emph{arithmeticity}: a complete classification of higher rank lattices. In the case where the semi-simple group is not almost simple, superrigidity reads as follows (see~\cite{Margulis}, page~2):

\begin{thm}[Margulis]
\label{thm_Margulis}%
Let $\Gamma$ be an irreducible lattice in $G = \prod_{\alpha\in A}\mathbf{G}_\alpha(k_\alpha)$, where $k_\alpha$ are local fields and $\mathbf{G}_\alpha$ are connected simply connected semi-simple $k_\alpha$\ti groups without $k_\alpha$\ti anisotropic factors and $|A|\geq 2$. Let $k$ be a local field, $\mathbf{H}$ a connected adjoint $k$\ti simple $k$\ti group and $\tau:\Gamma\to \mathbf{H}(k)$ a homomorphism with Zariski-dense unbounded image.

Then $\tau$ extends to a continuous homomorphism $\wt\tau: G\to \mathbf{H}(k)$.
\end{thm}

Our goal is to abandon completely the realm of algebraic (or Lie) groups and to establish a generalization of this theorem for uniform lattices in products. There will be no assumptions on the product group $G$. Instead of the algebraic group $\mathbf{H}$ for the target we shall consider general isometry groups of \cat0 spaces. An additional feature of our proof is that it is fully self-contained; this gives in particular a new, rather elementary, and purely geometric proof of Margulis' result for uniform lattices. For instance, the idiosyncrasies of positive characteristics (see Venkataramana~\cite{Venkataramana88}) vanish.

\smallskip

Here is our setting: \textbf{(i)}~Recall that a lattice $\Gamma$ in a product $G=G_1\times \cdots\times G_n$ of \textbf{arbitrary} locally compact groups is said \emph{irreducible} if the projection of $\Gamma$ to each factor $G_i$ is dense. In the classical semi-simple case of Theorem~\ref{thm_Margulis}, this follows from the stronger notion of algebraic irreducibility assumed therein.

\begin{rem*}
The irreducibility assumption is not a restriction whatsoever: One verifies that \emph{any} lattice $\Gamma<G$ is an irreducible lattice in the product $G^*<G$ of the closures $G^*_i<G_i$ of its projections to $G_i$. In particular, for any discrete cocompact subgroup $\Gamma<G$, the theorems below apply and extend the $\Gamma$\ti actions to $G^*$.
\end{rem*}

\textbf{(ii)}~We shall replace $\mathbf{H}(k)$ with an isometry group $\isom(X)$ of an \textbf{arbitrary} complete \cat0 space $X$; in the classical case, $X$ is the symmetric space of $\mathbf{H}(k)$ when $k$ is Archimedean and the associated Bruhat-Tits building otherwise. The Zariski-density is a necessary assumption for Theorem~\ref{thm_Margulis} as stated; before proposing two replacements for that assumption, we first state a result without any further assumption, which is possible for proper spaces upon passing to the geometric boundary:

\begin{thm}
\label{thm_boundary}%
Let $\Gamma$ be an irreducible uniform lattice in a product $G=G_1\times \cdots\times G_n$ of locally compact $\sigma$\ti compact groups. Let $\Gamma$ act by isometries on a proper \cat0 space $X$ without global fixed point.

Then there is a non-empty closed $\Gamma$\ti invariant set $\C\se\p X$ on which the $\Gamma$\ti action extends continuously to a $G$\ti action. Moreover this action factors through $G\to G_i$ for some $i=1,\ldots,n$.
\end{thm}

\noindent
(See the text for a more precise statement; \emph{e.g.} $\C$ arises as the boundary of a $G_i$\ti space lying in $X$.)

We turn now to extending homomorphisms $\tau : \Gamma \to H<\isom(X)$ as in Margulis' statement. We propose first, still for $X$ proper, the following substitute for Zariski-density in simple adjoint groups:

\begin{defi}
\label{defi_decomposable}%
A subgroup $L< \isom(X)$ is \emph{indecomposable} if for every non-empty $L$\ti invariant closed subset $\C\se \p X$, the closure $\bbb L$ acts faithfully on $\C$ and $\bbb L$ is closed in $\homeo(\C)$ (for the topology of uniform convergence).
\end{defi}

This always holds in the setting of Theorem~\ref{thm_Margulis}. Indecomposability turns out to be very natural, see Section~\ref{sec_decomposable}. In particular, a non-trivial indecomposable subgroup $L$ cannot fix any point at infinity. In negative curvature, indecomposability is essentially automatic upon passing to invariant subspaces.

\begin{cor}
\label{cor_superrigidity}%
Let $\Gamma$ be an irreducible uniform lattice in a product $G=G_1\times \cdots\times G_n$ of locally compact $\sigma$\ti compact groups. Let $H<\isom(X)$ be a closed subgroup, where $X$ is a proper \cat0 space, and let $\tau:\Gamma\to H$ be a homomorphism with indecomposable unbounded image.

Then $\tau$ extends to a continuous homomorphism $\wt\tau: G\to H$.
\end{cor}

This result immediately implies Margulis' Theorem~\ref{thm_Margulis} for uniform lattices.

Further, if we keep $H=\mathbf{H}(k)$, it shows similarly that for an irreducible uniform lattice in a \textbf{general} product group $G$, \itshape all completely reducible linear representations in finite-dimensional vector spaces over all local fields are completely determined by the continuous linear representations of $G$. \upshape Specialising in the other direction: Even when $G$ is an algebraic group, the above theorem yields a new family of superrigidity results.

\begin{rem*}
It is easy to verify (and inherent in the proofs) that the extended map $\wt\tau$ factors through some $G_i$. Similarly for Margulis' Theorem~\ref{thm_Margulis} and Theorem~\ref{thm_non-proper} below.
\end{rem*}

\medskip

Although the above results are set in the context of locally compact spaces, our proof involves in an essential way \textbf{infinite-dimensional \cat0 spaces}. Indeed, the overall strategy is to \emph{induce} the $\Gamma$\ti action to a $G$\ti action on a space of $L^2$~maps $G/\Gamma\to X$ and then to prove a splitting theorem for the latter. It turns out that we can prove superrigidity also when $X$ itself is infinite-dimensional; to this end, we propose our second substitute for Zariski-density in intrinsic geometric terms for \emph{any} \cat0 space $X$:

\begin{defi}
\label{defi_reduced}%
A subgroup $L< \isom(X)$ is \emph{reduced} if there is no unbounded closed convex subset $Y\subsetneqq X$ such that $g Y$ is at finite (Hausdorff) distance from $Y$ for all $g \in L$.
\end{defi}

\begin{thm}
\label{thm_non-proper}%
Let $\Gamma$ be an irreducible uniform lattice in a product $G=G_1\times \cdots\times G_n$ of locally compact $\sigma$\ti compact groups. Let $H<\isom(X)$ be a closed subgroup, where $X$ is any complete \cat0 space not isometric to a finite-dimensional Euclidean space~${\bf R}^d$. Let $\tau:\Gamma\to H$ be a homomorphism with reduced unbounded image.

Then $\tau$ extends to a continuous homomorphism $\wt\tau: G\to H$.
\end{thm}

\noindent
(\emph{Continuous} homomorphism and \emph{closed} subgroup are to be defined suitably when $X$ is not proper.)

\begin{rem*}
As stated, the above theorem does not hold for the very special case $X={\bf R}^d$ (see Section~\ref{sec_Rd}). Although Hilbert spaces, in particular ${\bf R}^d$, are special examples of complete \cat0 spaces, their linear structure allows a more detailed analysis; see Shalom~\cite{Shalom00}.
\end{rem*}

We have so far considered only uniform lattices. The non-uniform case raises some difficulties, which can however be overcome under certain assumptions:

\begin{thm}
\label{thm_non_uniform}%
Theorem~\ref{thm_boundary}, Corollary~\ref{cor_superrigidity}, Theorem~\ref{thm_non-proper}, and in fact all results of this paper hold more generally also for non-uniform lattices provided they are \emph{square-integrable} and \emph{weakly cocompact}. This is notably the case (i)~for all Kazhdan Kac-Moody lattices~\cite{Remy04} and (ii)~whenever $G$ is a connected semisimple Lie group.
\end{thm}

For a discussion of these concepts, more precise statements and proofs, see Appendix~B.

\subsection{Splitting} The geometry of infinite-dimensional \cat0 spaces is in some regards very different from their classical analogues; a first glimpse into the richness of phenomena arising there is given by the study of unitary representations and their cocycles -- Kazhdan's property~(T) is perhaps one of the first instances. A central issue is that the boundary $\p X$ does not reflect sufficiently the structure of the isometries of $X$; indeed $\p X$ may even be empty. Therefore, we introduce the following notion, which coincides simply with the existence of a fixed point in $\p X$ when $X$ is a proper \cat0 space:

\begin{defi}
\label{defi_eva}%
Let $G$ be a topological group with a continuous action by isometries on a metric space $X$. The $G$\ti action on $X$ is \emph{evanescent} if there is an unbounded set $T\se X$ such that for every compact set $Q\se G$ the set $\{d(g x, x): g\in Q, x\in T\}$ is bounded.
\end{defi}

We now can turn to our splitting theorem. Observe that there is no assumption whatsoever on the topology of $G$ or $X$:

\begin{thm}
\label{thm_splitting}%
Let $X$ be a complete \cat0 space with a continuous $G$\ti action by isometries, where $G=G_1\times \cdots\times G_n$ is any product of topological groups $G_i$.

\nobreak
Either the $G$\ti action is evanescent or there is a canonical non-empty closed convex $G$\ti invariant subspace $Z\se X$ which splits $G$\ti equivariantly isometrically as a product $Z_1\times\cdots\times Z_n$ of $G_i$\ti spaces $Z_i$.
\end{thm}

\noindent
(See Proposition~\ref{prop_Z2} for the importance of evanescence.)

\smallskip

The special case where $X$ is locally compact, though it is not representative of the above, is essentially known:

\begin{cor}
\label{cor_geometric}%
Let $X$ be a proper \cat0 space with a $G$\ti action by isometries, where $G=G_1\times \cdots\times G_n$ is any product of groups $G_i$.

\nobreak
Either there is a $G$\ti fixed point in $\p X$, or there is a non-empty closed convex $G$\ti invariant subspace $Z\se X$ which splits $G$\ti equivariantly isometrically as a product of $G_i$\ti spaces $Z_i$.\hfill\qedsymbol
\end{cor}

When $X$ is in fact a Riemannian manifold, this result is Schroeder's generalization~\cite{Schroeder} of the splitting theorems proved by Lawson-Yau~\cite{Lawson-Yau} and Gromoll-Wolf~\cite{Gromoll-Wolf} around~1970. The latter are of differential nature~-- and Schroeder's proof relies in an essential way on a Riemannian argument of Eberlein~\cite{Eberlein} using currents. However, there is a purely \cat0 statement in~\cite{Bridson-Haefliger} (II.6.21 and~II.6.25(3)); since the argument therein requires to extend geodesics indefinitely within every invariant subspace, it is assumed in that reference that $X$ has this property and that $G$ is cocompact or at least has \emph{full limit set}. Similar particular cases are obtained by Jost-Yau~\cite{Jost-Yau99} using harmonic maps. In our situation, it is impossible to assume anything on the space $Z$, but it turns out that no assumption is needed.

\subsection{Comments} One of the tools allowing us to deal with the infinite-dimensional spaces is a weakened topology \tc\ that we introduce; it is a common generalization of the weak topology in Hilbert spaces and the topology introduced in~\cite{Monod-Shalom1} for trees. Another more technical aspect is the behaviour of evanescence under induction (Theorem~\ref{thm_ind_eva} and Appendix~B).

\medskip

Other applications of our methods to rigidity theory will be exposed elsewhere; in particular, the special case of homomorphism to algebraic groups $\mathbf{H}(k)$ leads, following Margulis' ideas, to an \textbf{Arithmeticity \emph{vs.}\! Non-Linearity Alternative} for irreducible lattices in suitable product groups~\cite{MonodARITH}.

\medskip

\textbf{Related results.~}Margulis' Theorem~\ref{thm_Margulis} has been followed by numerous related results. We refer \emph{e.g.} to~\cite{BurgerICM},\cite{Gromov-Pansu} and references therein; to Zimmer's non-linear superrigidity~\cite{Zimmer84}; to the geometric superrigidity of Jost-Yau~\cite{Jost-Yau93b},\cite{Jost-Yau97},\cite{JostETH} and Mok-Siu-Yeung~\cite{Mok-Siu-Yeung}. Moreover, in the \cat{-1} setting, the flexibility of bounded cohomology allows for very general conclusions~\cite{Monod-Shalom1},\cite{Monod-ShalomCRAS}, including cocycle superrigidity in the spirit of Zimmer. By contrast, the present approach is more simple-minded: No boundary maps, no cohomology, no harmonic maps, and of course no theory of algebraic groups.

We also point out that an intermediate result in our strategy (Theorem~\ref{thm_superrigidity2}) provides a totally geodesic $\Gamma$\ti map in complete generality (from an auxiliary $G$\ti space), which is often enough to deduce geometric superrigidity.

\medskip

\textbf{On the exposition and use of previous work.~}We seek a self-contained presentation; in order to illustrate the new concepts that we introduce, we bring a number of examples or simple propositions which are not needed for the proof of the main results. We are inspired by various sources: Our first motivation is of course Margulis' work. The idea of inducing actions is classical in rigidity. Spaces of $L^2$ maps appear notably in~\cite{Korevaar-Schoen}. Our general strategy is analogous to Shalom's work for unitary representations and cohomology~\cite{Shalom00}. The proof of Theorem~\ref{thm_splitting} borrows some arguments from~\cite{Schroeder} and some others from~\cite[II.6]{Bridson-Haefliger}.

\medskip

\textbf{Acknowledgments.~}I am grateful to J.-F.~Quint for providing the elegant proof of Lemma~\ref{lemma_Zariski}; I thank D.~Fisher and T.~Gelander for helpful remarks on an earlier version. This work was partially supported by FNS grant~8220-067641 and NSF grant DMS~0204601.

\section{Informal Outline of the Reasoning}
%
Since our proof is accompanied with many geometric generalities, it might be useful to outline the general argument in order to highlight its simplicity:

\subsection{Superrigidity}
Consider a lattice $\Gamma<G = G_1\times \cdots\times G_n$ acting on a complete \cat0 space $X$.

\smallskip

We construct an associated \cat0 $G$\ti space $Y$ by considering the (right-) $\Gamma$\ti equivariant maps $f:G\to X$. The distance between $f,f'\in Y$ is $d^2(f,f') \defq \int_{G/\Gamma} d^2(f(g), f'(g)) \d g$ and the $G$\ti action is $(h f)(g) \defq f(h^{-1} g)$.

\smallskip

The motivation for this new ``induced'' \cat0 space $Y$ is this: On the one hand, there is a correspondence between properties of the $\Gamma$\ti action on $X$ and of the $G$\ti action on $Y$; to wit, non-evanescence is preserved. On the other hand, even though $Y$ is \emph{a priori} a more unwieldy space (typically not locally compact), the $G$\ti action on it is subject to the splitting theorem.

\smallskip

We obtain a splitting $Y\supseteq Z = Z_1 \times \cdots \times Z_n$ into $G_i$\ti spaces. If $f,f'\in Z_i$, then for $h\in G_j$ ($j\neq i$) we have a rectangle $\{f, f', h f, h f'\}$. It follows that $d(f(g), f'(g)) = d(f(h^{-1}g), f'(h^{-1}g))$ for almost all $g\in G$. If the projection of $\Gamma$ to $G_i$ is dense, this implies that $d(f(g), f'(g))$ is a.e. constant. We thus obtain isometric maps $Z_i\to X$ compatible with the $\Gamma$\ti actions. This is the main step in extending the $\Gamma$\ti action on $X$ or on $\p X$ to $G$.

\subsection{Splitting}
Consider a $G=G_1\times \cdots \times G_n$\ti action on a complete \cat0 space $Y$; we are motivated by the above ``induced'' spaces but work in complete generality.

\smallskip

Reduce to the case $n=2$. We show that either the $G$\ti action is evanescent or there is a minimal non-empty closed convex $G_1$\ti invariant subset $Z_1\se Y$. The main point here is a geometric analogue of the Banach-Alao\u{g}lu theorem: We introduce a weakened topology for which bounded closed convex sets are compact. A compactness argument then produces the minimal set.

\smallskip

The ``sandwich lemma'' implies as in~\cite{Schroeder} and~\cite[II.6]{Bridson-Haefliger} that the collection of all such sets $Z_1$ has a foliated structure which is preserved by the $G$\ti action. An addition to their arguments here is that even though we lack any further assumptions, this foliation is a \emph{global} isometric splitting. This follows by showing that the holonomy consists of Clifford translations, hence is trivial.

\section{Geometric Preliminaries}
\label{sec_quoi}%
\subsection{}
Our general background reference is~\cite{Bridson-Haefliger} (see also~\cite{Ballmann},\cite{JostETH}). Let $X$ be a metric space with metric $d$. A map $\sigma:I\to X$, where $I\se{\bf R}$ is any interval, is \emph{geodesic} if it is isometric. The space $X$ is \emph{geodesic} if every pair of points is joined by some geodesic segment, and $A\se X$ is \emph{convex} if it contains any geodesic segment joining any two of its points.  For any subsets $A,B$ of a metric space $X$, we denote by $[A]$ the closed convex hull of $A$ and $[A,B]\defq[A\cup B]$. We write $\cball(x,r)$ for the closed ball of radius $r$ around $x$; if $X$ is a geodesic space and $r\neq 0$, then this is the closure of the open ball $\ball(x,r)$.

\subsection{}
\label{sec_cat}%
A \emph{CAT(0) space} is a geodesic space $X$ such that for any triangle $x,c,c'\in X$ the midpoint $m$ of any geodesic between $c,c'\in X$ satisfies the \emph{courbure n\'egative} inequality of Bruhat-Tits~\cite[II.1.9]{Bridson-Haefliger}:
\beq
\label{eq_CN}%
2 d^2(m,x)\ \leq\ d^2(c',x) + d^2(c,x)  - \frac{1}{2} d^2(c',c).
\eeq
Equivalently, the distance between any points in the sides of any geodesic triangle are bounded by the corresponding distances in Euclidean triangles. In particular, geodesics are unique: the segment from $c$ to $c'$ is $[c,c']$. Examples include all symmetric spaces and Bruhat-Tits buildings; all simply connected manifolds of non-positive sectional curvature; Hilbert spaces; simply connected Euclidean or hyperbolic simplicial complexes satisfying certain local link conditions~\cite[II.5.4]{Bridson-Haefliger}.

Here is an example of application of the \cat0 inequalities. Let $x,x',y,y'\in X$. For $0<\epsi<1$ let $x_\epsi, x'_\epsi$ be the points of $[x,x']$ at distance $\epsi d(x,x')$ of $x$, respectively of $x'$. Then
\beq
\label{eq_Resh}%
d^2(x_\epsi,y) + d^2(x'_\epsi,y')\ \leq\ d^2(x,y) + d^2(x',y') + 2\epsi d(x,x') \big( d(y,y') - (1-\epsi) d(x,x')\big).
\eeq
This type of inequalities was first proved by Reshetnyak~\cite{Reshetnyak68} (we refer to~\cite[2.1.3]{Korevaar-Schoen} for a derivation of~(\ref{eq_Resh}) from~(\ref{eq_CN})).

\subsection{}
The \emph{circumradius} of a bounded set $A\se X$ is the infimum $\ro$ of all $r>0$ such that there is $x\in X$ with $A\se \cball(x,r)$. If $X$ is \cat0 and complete, then this infimum is achieved and (for $A\neq\varnothing$) there is a unique point $c\in X$ such that $A\se \cball(c,\ro)$; this point is called the \emph{circumcentre} of $A$. If in addition $A$ is convex and closed, then $c\in A$.

\begin{lemma}
\label{lemma_circum}%
Let $X$ be a complete \cat0 space, $E\se E'\se X$ two non-empty bounded closed convex sets, $c,c'$ the corresponding circumcentres and $\ro,\ro'$ the circumradii.

Then $d(c',c)\leq \sqrt{2}\sqrt{{\ro'}^2-\ro^2}$.
\end{lemma}

\begin{proof}
If $c=c'$, there is nothing to prove. Otherwise, the midpoint $m$ of $[c',c]$ is not the circumcentre of $E$ and thus there is $x\in E$ with $d(m,x)>\ro$. Replacing this in~(\ref{eq_CN}) gives $2 \ro^2 + d^2(c',c)/2 < {\ro'}^2 + \ro^2$, as required.
\end{proof}

\subsection{}
\label{sec_sandwich}%
If $A\se X$ a non-empty closed convex subset of the complete \cat0 space $X$, then there is a \emph{nearest point projection} map $p_A: X\to A$ which does not increase distances.

Another consequence of~(\ref{eq_CN}) is the Sandwich Lemma (see~\cite[II.2.12]{Bridson-Haefliger}): Assume that $C,C'\se X$ are two non-empty closed convex subsets such that the function $d(x,C)$ is constant on $x\in C'$ and likewise $d(x,C')$ on $C$; denote this (common) constant by $d_0$. Then there is a canonical isometry
$$\fhi:\ C\times [0,d_0] \xrightarrow{\ \cong\ }[C,C']\se X$$
such that
$$\fhi|_{C\times\{0\}} = \id_C\ \mbox{ and }\ \fhi|_{C\times\{d_0\}} = p_{C'}|_C.$$
In the particular case where $C,C'$ are geodesic segments, we say that they \emph{determine a Euclidean rectangle}.

\medskip

Similar arguments apply to the following setting: An isometry $g:X\to X$ such that $D\defq d(gx, x)$ is independent of $x$ is called a \emph{Clifford translation}~\cite[II.6.14]{Bridson-Haefliger}. If $D>0$, then there is an isometry $X\cong X' \times {\bf R}$ intertwining $g$ with the translation by $D$ along ${\bf R}$. More generally, any complete \cat0 space $X$ splits canonically as $X\cong X'\times V$, where $V$ is a (possibly trivial or finite-dimensional) Hilbert space and $X'$ does not admit any Clifford translation. Moreover, $\isom(X)$ preserves this splitting. For all this, see~\cite[II.6.15]{Bridson-Haefliger}.

\subsection{}
\label{sec_group_action}%
Let $G$ be a topological group acting on a metric space $X$ by isometries. The following are equivalent:

(i)~For all $x\in X$ the map $G\to X$, $g\mapsto gx$ is continuous at $e\in G$.

(ii)~The action map $G\times X\to X$ is continuous.

\noindent
When this happens, we say that the action is \emph{continuous}. When $X$ is proper, we always endow its isometry group $\isom(X)$ with the compact-open topology; in that case, $\isom(X)$ is a locally compact second countable topological group and the above conditions are equivalent to the continuity of the homomorphism $G\to\isom(X)$. For $X$ general, we do not topologize $\isom(X)$ but still call a homomorphism $G\to \isom(X)$ continuous when the $G$\ti action is so, and a subgroup $H<\isom(X)$ is said closed if its orbits in $X$ are so. A standard argument implies:

\begin{lemma}
\label{lemma_cont_Borel}%
Let $G$ be a locally compact second countable group with an action by isometries on a complete separable metric space $X$. If for all $x\in X$ the map $G\to X$, $g\mapsto g x$ is measurable, then the action is continuous.\hfill\qedsymbol
\end{lemma}

We call a subset $L\se G$ \emph{bounded} if for some (or equivalently any) $x\in X$ the set $Lx$ is bounded in $X$; when $X$ is proper, this coincides with the usual definition in which a subset of a locally compact group $G$ is bounded when it has compact closure.

\subsection{}
\label{sec_boundary}%
The \emph{boundary at infinity} $\p X $ of a complete \cat0 space $X$ can be defined as the set of equivalence classes of geodesic rays in $X$, where two rays are equivalent (\emph{asymptotic}) if they remain at bounded distance from each other. The \emph{bordification} $\bbb X \defq X\sqcup \p X$ can be identified with the inverse limit of the closed balls $\cball(x,r)$ as $r\to\infty$ under the maps $p_{\cball(x,r)}$, wherein $x\in X$ is any basepoint; the resulting topology is called the \emph{cone topology}~\cite[II.8]{Bridson-Haefliger}. When $X$ is proper, $\bbb X$ is a second countable compact space (hence metrizable) and the action map $\isom(X)\times \bbb X\to \bbb X$ is continuous.

\subsection{}
Let $X$ be a metric space. Unless otherwise stated, every topological epithet will always refer to the topology $\T$ induced by the metric $d$. One can define a weaker topology \tw\ by letting \tw\ be the weakest topology on $X$ such that for all $x,y\in X$ the map $z\mapsto d(x,z) - d(y,z)$ is continuous. This topology is always Hausdorff; we shall however be more interested in the following:

\begin{defi}
Let $X$ be a metric space. We define the topology \tc\ to be the weakest topology on $X$ for which all $\T$\ti closed convex sets are \tc\ti closed.
\end{defi}

Here is the main property of \tc.

\begin{thm}\addtocounter{footnote}{1}
\label{thm_compact}%
Let $X$ be a complete \cat0 space and $C\se X$ a bounded closed convex subset. Then $C$ is compact\footnote{we follow the common usage to define compactness with the Borel-Lebesgue axiom regardless of separation; this is called \emph{quasi-compact} by Bourbaki's collaborators~\cite[I \S9 \no1]{BourbakiTG12}. On the other hand, our locally compact groups are always assumed Hausdorff.} for the topology \tc.
\end{thm}

This is a common generalization of the two quite different cases of Hilbert spaces in the weak topology and of trees in the topology considered in~\cite{Monod-Shalom1}; see the examples below.

\begin{proof}[Proof of the theorem]
We need to prove that for any family $\F$ of \tc\ti closed sets $F\se X$ such that $\{F\cap C: F\in \F\}$ has the finite intersection property, the intersection of all $F\cap C$ is non-empty. By Alexander's sub-base theorem, it is enough to consider a family $\F$ consisting of closed convex sets $F$. We may assume $F\se C$ upon replacing each $F$ by $F\cap C$. Let $A$ be the set of non-empty finite subsets $\A\se \F$, ordered by inclusion. To any $\A\in A$ we associate the circumradius $\ro_\A$ of $\bigcap \A$ and its circumcentre $c_\A\in \bigcap \A$. Notice that $\{\ro_\A\}_{\A\in A}$ is a non-increasing net since $\A\se\B$ implies
$$\bigcap \B \se \bigcap \A \se \cball(c_\A, \ro_\A)$$
and thus $\ro_\B\leq \ro_\A$. On the other hand this net is non-negative, therefore converges and thus is a Cauchy net. Applying Lemma~\ref{lemma_circum} to the sets $\bigcap \B \se \bigcap \A$ we deduce that $\{c_\A\}_{\A\in A}$ is a Cauchy net and hence converges to a limit $c\in X$. For every $\A\in A$, all points $c_\B$ with $\A\se\B$ belong to $\bigcap \A$; therefore the limit $c$ is still in $\bigcap \A$. It follows that $c$ is in
$$\bigcap_{\A\in A}\bigcap \A = \bigcap_{F\in \F} F,$$
as was to be shown.
\end{proof}

\begin{rem}
\label{rem_non-expanding}%
Lemma~\ref{lemma_circum} can also be used to explain the following (probably well-known) fact: Let $F$ be a non-expanding map of a complete \cat0 space $X$. If $F$ has bounded orbits, then its (closed convex) fixed set is non-empty. Indeed, pick any $x\in X$; let $C_n \defq [\{F^k(x): k\geq n\}]$ and let $c_n$ be its circumcentre, $\ro_n$ its circumradius. Since $\{C_n\}$ is decreasing, Lemma~\ref{lemma_circum} implies that $\{c_n\}$ has a limit $c$. Since $C_{n+1}\se [F(C_n)]$, we have $C_{n+1} \se \cball(F(c_n), \ro_n)$. Applying again Lemma~\ref{lemma_circum}, it follows that $F(c_n)$ converges to $c$, which is thus fixed.
\end{rem}

\subsection{}
\label{sec_ill_Tc}%
This Section~\ref{sec_ill_Tc}  serves only to illustrate \tc.

\begin{lemma}
\label{lemma_top_incl}%
For any complete \cat0 space $X$ we have $\tc\se\tw\se \T$.
\end{lemma}

\begin{proof}
We only need to show $\tc\se\tw$. Let $C\subsetneqq X$ be closed convex. The set
$$C' \defq \bigcap_{x\notin C}\big\{z : d(x,z)-d(p_C(x),z)\geq 0\big\}$$
is \tw\ti closed. We claim that $C'=C$. On the one hand, $C\se C'$ because for all $z\in C$ we have $d(p_C(x),z)=d(p_C(x), p_C(z)) \leq d(x,z)$ regardless of $x$. On the other hand, if $z\notin C$ then $z$ is not in $\{z : d(x,z)-d(p_C(x),z)\geq 0\}$ for $x=z$.
\end{proof}

\begin{lemma}
\label{lemma_top_coincide}%
If $K\se X$ is $\T$\ti compact, then the restriction of $\T$ and $\tc$ to $K$ coincide.
\end{lemma}

\begin{proof}
It is enough to show that for all $x\in K$ and all $r>0$ there is a $\tc$\ti neighbourhood $U$ of $x$ such that $U\cap K \se \ball(x,r)$. That ball, together with all $\ball(y,r/2)$ for every $y\in K$ with $d(x,y) \geq r$, covers $K$. Therefore there is a finite set $F\se K$ of points $y$ with $d(x,y)\geq r$ such that $\bigcup_{y\in F} \ball(y,r/2)$ covers $K \setminus \ball(x,r)$. The set $U \defq X \setminus \bigcup_{y\in F} \cball(y,r/2)$ has the sought properties.
\end{proof}

The topology \tc\ is familiar in some particular cases:

\begin{exo}
If $X$ is a Hilbert space, then \tc\ and \tw\ both coincide with the weak topology.
\end{exo}

\begin{exo}
\label{exo_hyp}%
If $X$ is the standard infinite-dimensional (separable) real hyperbolic space $\mathbf{O}(1,\infty)/\mathbf{O}(\infty)$ (see \emph{e.g.}~\cite{Burger-Iozzi-Monod}), then \tw\ and \tc\ coincide; moreover, they are induced from the weak topology if we realize $X$ with the ball model in a Hilbert space.
\end{exo}

\begin{exo}
\label{exo_trees}%
If $X$ is a simplicial tree, then \tw\ and \tc\ coincide; moreover, they coincide with the weak topology $\sigma$ on trees introduced in~\cite{Monod-Shalom1}.
\end{exo}

However, the topology \tc\ is not all that straightforward in general. It seems not to be Hausdorff for higher rank Bruhat-Tits buildings or symmetric spaces, and it is unclear to us what happens for \emph{complex} hyperbolic spaces (of finite or infinite dimension).

Moreover, whilst in a Hilbert space all weakly compact sets are bounded, this is not so even in the simplest examples of \cat0 spaces:

\begin{exo}
Let $X$ be a simplicial tree consisting of countably many rays of finite but unbounded length, all meeting at one vertex. Then the space $X$ is \tc\ti compact even though unbounded. Notice in addition that $\p X=\varnothing$.
\end{exo}

Let $X$ be any complete \cat0 space. We extend the topology \tc\ to $\bbb X$ by declaring that for any $\T$\ti closed convex set $C\se X$ the (usual) closure $\bbb C$ of $C$ in $\bbb X$ is \tc\ti closed.

\begin{rem}
The (compact) topology that \tc\ determines on $\bbb X$ through the realization of $\bbb X$ as inverse limit of closed balls is in general coarser than \tc, even when restricted to $X$: Already for Hilbert spaces, one obtains the weaker \emph{bounded weak topology}. These topologies coincide however when $X$ is a tree or a real hyperbolic space as in Example~\ref{exo_hyp}.
\end{rem}

Whilst $\bbb X$ is not \tc\ti compact when $X$ is \emph{e.g.} an infinite-dimensional Hilbert space, we have:

\begin{prop}
If the complete \cat0 space $X$ is Gromov-hyperbolic, then $\bbb X$ is \tc\ti compact.
\end{prop}

\begin{proof}
It suffices to show that for any nested family $\C$ of non-empty closed convex sets $C\se X$ the intersection $\bigcap_{C\in\C} \bbb C$ is non-empty. Fix $x\in X$. If $\{p_C(x) : C\in\C\}$ is bounded, we are done by Theorem~\ref{thm_compact}. Otherwise, we claim that $\bigcap_{C\in\C} \bbb C$ is a single point of $\p X$. Indeed, in view of the sequential model of $\p X$~\cite[p.~120]{Ghys-Harpe}, it is enough to show that for any choices $c\in C$, $c'\in C'$, the Gromov product $(c|c')_x$ tends to infinity along $C,C'\in\C$. This follows since $(c|c')_x$ is comparable to $d([c,c'],x)$, which is bounded below by $d(p_C(x), x)$ or $d(p_{C'}(x), x)$.
\end{proof}

\subsection{}
\label{sec_eva}%
We will now analyse \emph{evanescence} (Definition~\ref{defi_eva}); let $G$ be a topological group with a continuous action by isometries on a metric space $X$.

\begin{defi}
\label{defi_eva_set}%
Let $Q$ be a subset of $G$. Then a subset $T\se X$ such that $\{d(g x, x): g\in Q, x\in T\}$ is bounded will be said \emph{$Q$\ti evanescent}.
\end{defi}

Thus, by Definition~\ref{defi_eva}, the action is evanescent if and only if there is an unbounded set $T\se X$ which is $Q$\ti evanescent for every compact set $Q\se G$; we then call $T$ itself \emph{evanescent}.

\begin{lemma}
\label{lemma_eva_lin}%
Suppose that $X$ is a \cat0 space, let $Q$ be compact in $G$ and $x_0\in X$. Then there exists no unbounded $Q$\ti evanescent set if and only if there is $\lambda>0$ and $d_0\geq 0$ such that
$$\sup_{g\in Q} d(gx, x) \geq \lambda d(x,x_0)-d_0\kern 1cm\forall\,x\in X.$$
\end{lemma}

\begin{proof}
Sufficiency is obvious. Suppose conversely that the condition fails. Then for every $n\geq 1$ there is $y_n\in X$ such that $d(g y_n, y_n) < d(y_n, x_0)/n! - n^2$ for all $g\in Q$. Let $x_n$ be the point at distance $n$ from $x_0$ on $[x_0, y_n]$. A comparison argument shows that
$$\limsup_{n\to\infty}\sup_{g\in Q} d(g x_n, x_n) \leq \sup_{g\in Q} d(g x_0, x_0)$$
and thus $\{x_n\}$ is a $Q$\ti evanescent sequence.
\end{proof}

\begin{rem}
\label{rem_ultra}%
It follows from this lemma that the $G$\ti action on $X$ is evanescent if and only if every compactly generated subgroup of $G$ has a non-trivial fixed point for its diagonal action on some (or equivalently any) \emph{asymptotic cone} of $X$ along a free ultrafiltre (the base-point of an asymptotic cone is a trivial fixed point). We shall not use this characterization.
\end{rem}

\begin{prop}
\label{prop_eva_proper}%
Suppose that $X$ is a complete \cat0 space.

\nobreak\noindent
(i)~If there exists a $G$\ti fixed point in $\p X$ then the $G$\ti action on $X$ is evanescent.

\nobreak\noindent
(ii)~The converse holds if $X$ is proper.
\end{prop}

\begin{proof}
If there exists a point $\xi\in\p X$ fixed by $G$, then any ray pointing to $\xi$ is an evanescent set. In case~(ii), $\bbb X=X\sqcup \p X$ is compact and it follows that any unbounded evanescent sequence has a subsequence converging to some $\xi\in\p X$; the definition of the cone topology on $\bbb X$ shows that $\xi$ is $G$\ti fixed.
\end{proof}

Another natural definition is as follows:

\begin{defi}
\label{defi_w_eva}%
A continuous $G$\ti action on a metric space $X$ is \emph{weakly evanescent} if for every compact $Q\se G$ there is an unbounded $Q$\ti evanescent set in $X$.
\end{defi}

We point out that for both Definitions~\ref{defi_eva} and~\ref{defi_w_eva} it is enough to consider unbounded evanescent \emph{sequences}. Whilst evanescence is in general stronger than weak evanescence, we have:

\begin{prop}
\label{prop_ceva}%
Suppose that $X$ is a complete \cat0 space. Then weak evanescence implies evanescence if either (i)~$X$ is proper; or~(ii) if there is a countable cofinal chain in the directed set of compact subsets of $G$.
\end{prop}

\begin{rem}
\label{rem_ceva}%
The assumption in~(ii) is satisfied if $G$ is locally compact $\sigma$\ti compact (since $G$ is then a countable union of compact sets of non-empty interior); but it also holds in other cases (\emph{e.g.} if $G$ is the additive group of a dual Banach space endowed with the \weak topology).
\end{rem}

\begin{proof}[Proof of the proposition]
For~(i), observe that the family $\big\{(\p X)^Q : Q \mbox{ compact in } G\big\}$ has the finite intersection property so that we can use compactness of $\p X$ and Proposition~\ref{prop_eva_proper}. In case~(ii), one can choose an increasing sequence of compact sets $Q_n\se G$ such that any compact set in contained in $Q_n$ for $n$ big enough. Let $T_n$ be unbounded and $Q_n$\ti evanescent and let $K_n\geq 1$ be a constant bounding the numbers $d(g x,x)$ of Definition~\ref{defi_eva_set}. Fix any point $x_0\in X$ and define a diagonal sequence $\{x_n\}$ by picking first $y_n\in T_n$ with $d(y_n, x_0)\geq (n K_n)^2$ and then letting $x_n$ be the point at distance $n$ of $x$ on $[x,y_n]$. A comparison argument as in the proof of Lemma~\ref{lemma_eva_lin} shows that $\{x_n\}$ is evanescent.
\end{proof}

Evanescence behaves in a simple way with respect to direct products:

\begin{prop}
\label{prop_eva_prod}%
Let $G=G_1\times G_2$ be a product of topological groups and suppose that $X=X_1\times X_2$ is a product of two unbounded \cat0 spaces $X_i$ with continuous $G_i$\ti action. Endow $X$ with the product $G$\ti action.

\nobreak\noindent
(i)~Both $G_i$\ti actions on $X$ are evanescent.

\nobreak\noindent
(ii)~The $G$\ti action on $X$ is (weakly) evanescent if and only if at least one of the $G_i$\ti actions on $X_i$ is (weakly) evanescent.
\end{prop}

Observe that~(i) stands in contrast to Shalom's notion of (non-)uniformity~\cite{Shalom00}.

\begin{proof}
For the first point, fix any $x\in X_1$ and let $T=\{x\}\times X_2$. This set is evanescent for $G_1$, since for any compact $Q\se G_1$ the set $Q x$ is bounded and $G_1$ acts trivially on $X_2$; likewise for $G_2$. A similar argument shows that if, say, the $G_1$\ti action on $X_1$ is (weakly) evanescent, then the $G$\ti action on $X$ is so too. As for the converse, assume $\{x_n\}$ is an unbounded $Q$\ti evanescent sequence in $X$ for some compact $Q\se G$. Fix $a_i\in X_i$ and set $x'_n = p_{X_1\times\{a_2\}}(x_n)$, $x''_n= p_{\{a_1\}\times X_2}(x_2)$. Since $p_{X_1\times\{a_2\}}$ is $G_1$\ti equivariant and does not increase distances (and similarly for $p_{\{a_1\}\times X_2}$), it is enough to show that either $\{x'_n\}$ or $\{x''_n\}$ is unbounded. But otherwise $\{x_n\}$ would itself be bounded.
\end{proof}

\subsection{}
\label{sec_Kazhdan}%
We recall that following Kazhdan a unitary or orthogonal representation $\pi$ is said \emph{almost to have non-zero invariant vectors} if for every compact $Q\se G$ and every $\epsi>0$ there is a unit vector $v$ such that
\beq
\label{eq_Kazhdan}%
\sup_{g\in Q}\|\pi(g)v - v\| \leq \epsi.
\eeq

\begin{prop}
\label{prop_Kazhdan}%
Suppose that $X$ is a Hilbert space. Then the $G$\ti action on $X$ is weakly evanescent if and only if the associated orthogonal representation almost has non-zero invariant vectors.
\end{prop}

\begin{proof}
A straightforward verification using \emph{e.g.} sequences of the form $x_n=nv_n$, where $v_n$ satisfies~(\ref{eq_Kazhdan}) for $\epsi= 1/n$.
\end{proof}

\begin{rem}
\label{rem_asympt_inv}%
In particular, one has the following well known characterization when $G$ is as in Remark~\ref{rem_ceva}: $\pi$ almost has non-zero invariant vectors if and only if there is an \emph{asymptotically invariant sequence}, \emph{i.e.} a sequence $\{v_n\}$ with $\|v_n\|=1$ and $\|g v_n - v_n\|\to0$ uniformly for $g$ on compact sets.
\end{rem}

\section{A General Splitting Theorem}
\label{sec_splitting}%
This section addresses Theorem~\ref{thm_splitting}.

\subsection{}
\label{sec_Z2}%
Consider the very simplest product group: $G={\bf Z}\times {\bf Z}$, and the very simplest non-proper \cat0 space: a Hilbert space. Then, already, evanescence is the right replacement for fixed points at infinity in the sense that Theorem~\ref{thm_splitting} fails if we formulate it in terms of $\p X$. We give a more general counter-example:

\begin{prop}
\label{prop_Z2}%
Let $G_1, G_2$ be any two locally compact separable groups without Kazhdan's property~(T). Then there is an action of $G=G_1\times G_2$ by isometries on a Hilbert space $X$ such that:

\smallskip\noindent
(i)~The $G$\ti action on $\p X$ has no fixed point.

\noindent
(ii)~No non-empty closed convex $G$\ti invariant subspace $Z\se X$ splits as product of $G_i$\ti spaces.
\end{prop}

\begin{proof}
By Theorem~1 in~\cite{Bekka-Valette}, each $G_i$ admits an orthogonal representation $\pi_i$ such that $\pi_i\otimes \pi_i$ almost has non-zero invariant vectors but $\pi_i$ is \emph{weakly mixing}, that is, has no finite-dimensional subrepresentation. Let $\sigma$ be the $G$\ti representation $\pi_1\otimes \pi_1\otimes \pi_2\otimes \pi_2$ and $X$ the associated Hilbert space. Notice that $\pi_i\otimes\pi_i$ itself is weakly mixing; this implies $\sigma^{G_i}=0$, see~\cite{Bergelson-Rosenblatt}. On the other hand, one deduces immediately with the definition~(\ref{eq_Kazhdan}) that $\sigma$ almost has non-zero invariant vectors. As a well-known consequence of the closed graph theorem (observed by Guichardet~\cite[Th\'eor\`eme~1]{GuichardetII}), it follows that there is a non-trivial $\sigma$\ti cocycle $b:G\to X$ in the closure of the space of coboundaries. We endow $X$ with the corresponding continuous $G$\ti action by (affine) isometries. For~(i), observe that a fixed point in $\p X$ would give a non-zero $G$\ti fixed vector for $\sigma$. For~(ii), assume for a contradiction that $Z\se X$ splits as $Z_1\times Z_2$. For every $x,y\in Z_i$, the vector $x-y$ is fixed by $\sigma(G_j)$ for $j\neq i$, which in view of $\sigma^{G_j}=0$ shows that both $Z_i$ are reduced to a single point. This point being fixed by the affine $G$\ti action, $b$ is trivial, contradicting the assumption.
\end{proof}

Notice that the action here is weakly evanescent by Proposition~\ref{prop_Kazhdan}, hence evanescent by Proposition~\ref{prop_ceva}. Also, the action is \emph{non-uniform} in Shalom's sense~\cite{Shalom00}, and indeed the above situation stands in contrast to the results of~\cite{Shalom00} for uniform actions.

\subsection{}
We undertake now the proof of Theorem~\ref{thm_splitting}. First observe that it is enough to consider the case $n=2$. Indeed, the case $n=1$ is tautological; further, assume that for $n\geq 3$ we have a subspace $Z$ splitting equivariantly as $Z=Z_1\times Z'_1$, where $Z_1$ is a $G_1$\ti space and $Z'_1$ a $G'_1=G_2\times\cdots\times G_n$\ti space. In order to apply induction, we just need to observe that the $G'_1$ action on $Z'_1$ cannot be evanescent unless the $G$\ti action on $X$ is so, see Proposition~\ref{prop_eva_prod}.

We can from now on assume that the $G$\ti action on $X$ is not evanescent. Theorem~\ref{thm_compact} allows us to get started with the following:

\begin{prop}
\label{prop_existence}%
Let $X$ be a complete \cat0 space with a continuous $G$\ti action by isometries, where $G=G_1\times G_2$ is any product of two topological groups.

\nobreak
If the $G$\ti action is not evanescent then there is a minimal non-empty closed convex $G_1$\ti invariant set in $X$.
\end{prop}

\begin{proof}
Choose any point $x\in X$ and let $\C$ be the set of non-empty closed convex $G_1$\ti invariant subsets of $[G_1 x]$, ordered by reverse inclusion. By Hausdorff's maximal principle, there is a (non-empty) maximal totally ordered subset $\D\se \C$. If the intersection $\bigcap\D$ is non-empty, we are done; we assume for a contradiction that it is empty and will show that the $G$\ti action on $X$ is evanescent. The net $d(D,x)$ indexed by $D\in\D$ is non-decreasing; we shall prove that it converges to infinity.

Indeed, if not, it would have a limit $d\in {\bf R}$ and the family $\big\{D\cap \cball(x,d): D\in \D\big\}$ would be a nested collection of non-empty closed convex bounded sets with empty intersection. This contradicts \tc\ti compactness established in Theorem~\ref{thm_compact}. It follows that the set $\{p_D(x):D\in\D\}$ is unbounded. It is now enough to show that for any compact sets $Q_i\se G_i$ the set
\beq
\label{eq_prod_ev}%
\big\{ d(g_1g_2 p_D(x), p_D(x)) : g_i\in Q_i, D\in\D\big\}
\eeq
is bounded. On the one hand, $p_D$ is $G_1$\ti equivariant and thus 
$$d(g_1 p_D(x), p_D(x)) = d(p_D(g_1 x), p_D(x)) \leq d(g_1 x, x)$$
implies that the family $d(g_1 p_D(x), p_D(x))$ is bounded. On the other hand, for any $g_2\in G_2$ the function $y\mapsto d(g_2 y, y)$ is constant on $G_1$\ti orbits, so that by convexity and continuity of the metric it is bounded by $d(g_2 x, x)$ on $y\in[G_1 x]$. It follows now from
$$d(g_1g_2 p_D(x), p_D(x)) \leq d(g_1 p_D(x), p_D(x)) + d(g_2 p_D(x), p_D(x))$$
that the collection~(\ref{eq_prod_ev}) is bounded. 
\end{proof}

\begin{rem}
\label{rem_existence}%
In particular, any non-evanescent $G$\ti action admits a minimal non-empty closed convex $G$\ti invariant set and actually any nested family of non-empty closed convex $G$\ti invariant sets has non-empty intersection. See also Remarks~\ref{rems_exi_bis}.
\end{rem}

\subsection{}
We now know that the set $Z_2$ of minimal non-empty closed convex $G_1$\ti invariant sets in $X$ is non-empty. We shall use repeatedly the following obvious

\begin{lemma}
\label{lemma_const}%
If $C\se X$ is a minimal non-empty closed convex $G_1$\ti invariant set, then any convex $G_1$\ti invariant continuous (or lower semi-continuous) function on $C$ is constant.
\end{lemma}

\begin{proof}
If $f:C\to{\bf R}$ were to assume two distinct values $s<t$, then $\big\{x\in C : f(x)\leq s\big\}$ would be a strictly smaller non-empty closed convex $G_1$\ti invariant set.
\end{proof}

Fix $Z_1\in Z_2$ and let $Z \defq \bigcup Z_2\se X$. Observe that each $Z_i$ has a natural $G_i$\ti action. We can now consider the following setup, borrowed from~\cite[pp.~239--241]{Bridson-Haefliger}: For every $C,C'\in Z_2$, the function $d(x,C)$ is constant on $C'$ by Lemma~\ref{lemma_const} since it is $G_1$\ti invariant; likewise with $C,C'$ interchanged. The Sandwich Lemma (Section~\ref{sec_sandwich}) yields a canonical isometry
\beq
\label{eq_fhi}%
\fhi:\ C\times [0,d(C,C')] \xrightarrow{\ \cong\ }[C,C']\se X
\eeq
such that
\beq
\label{eq_fhi_prop}%
\fhi|_{C\times\{0\}} = \id_C\ \mbox{ and }\ \fhi|_{C\times\{d(C,C')\}} = p_{C'}|_C.
\eeq
In particular, the distance $d(C,C')$ defines indeed a metric on the set $Z_2$ and this metric is geodesic. Furthermore we have $Z=\bigsqcup_{C\in Z_2}C$ and hence obtain a well-defined bijection $\alpha:\ Z\lra Z_1\times Z_2$ by setting $\alpha(x)=(p_{Z_1}(x), C_x)$, wherein $C_x$ is the unique element of $Z_2$ containing $x$.

\subsection{}
At this point, the main remaining steps are to show that $\alpha$ is actually an isometry and that it intertwines the $G$\ti actions; we need the following key fact:

\begin{prop}
\label{prop_coherent}%
For all $C_1, C_2, C_3\in Z_2$ we have $p_{C_1}\circ p_{C_3}\circ p_{C_2}|_{C_1} = \id_{C_1}$.
\end{prop}

\begin{proof}
First we point out that if all three sets $C_i$ were just geodesic lines, then this proposition is a well known general fact holding for any three parallel lines in any metric space, see~\cite[II.2.15]{Bridson-Haefliger}. However, in our case, it is not even necessary that $C_i$ should contain any line. Thus, denote by $\teta:C_1\to C_1$ the above map $p_{C_1}\circ p_{C_3}\circ p_{C_2}|_{C_1}$. The properties of the isometry $\fhi$ in~(\ref{eq_fhi}) imply that $p_{C_j}|_{C_i}$ is an isometry for all $i,j$; moreover, this isometry $C_i\to C_j$ is $G_1$\ti equivariant because $C_j$ is $G_1$\ti invariant. Therefore, $\teta$ is a $G_1$\ti equivariant isometry. It follows that the function $x\mapsto d(\teta(x), x)$ is a convex continuous $G_1$\ti invariant function of $x\in C_1$; Lemma~\ref{lemma_const} implies that it is constant. Thus $\teta$ is a Clifford translation of $C_1$. We need to show that $\teta$ is trivial. But a non-trivial Clifford translation preserves the image of a geodesic line $\sigma:{\bf R}\to C_1$; indeed, recall that in fact in that case $C_1$ would split for $\teta$ (Section~\ref{sec_sandwich}; see~\cite[II.6.15]{Bridson-Haefliger}~-- though we will not need this. We may now apply the general fact mentioned earlier (\cite[II.2.15]{Bridson-Haefliger}) to the three lines $\sigma$, $p_{C_2}\circ\sigma$ and $p_{C_3}\circ \sigma$. We deduce that $\teta$ translates $\sigma({\bf R})$ trivially. Therefore, the constant $d(\teta(x), x)$ vanishes and $\teta=\id_{C_1}$ as was to be shown.
\end{proof}

\subsection{}
It follows now that $\alpha$ is isometric; for completeness (and because it contains a misprint), we give the calculation of~\cite[p.~241]{Bridson-Haefliger}: Let $x,x'\in Z$; using twice that~(\ref{eq_fhi_prop}) defines an isometry to the Cartesian product $C\times [0,d(C,C')]$ for all $C,C'\in Z_2$, we have
$$d^2(x,x') = d^2(x, p_{C_x}(x')) + d^2(C_x, C_{x'}) = d^2(p_{Z_1}(x), p_{Z_1}\circ p_{C_x}(x')) + d^2(C_x, C_{x'}).$$
Applying now Proposition~\ref{prop_coherent} to $C_x, C_{x'}$ and $Z_1$ we deduce $p_{Z_1}\circ p_{C_x}(x') = p_{Z_1}(x')$, so that $d^2(x,x') = d^2(p_{Z_1}(x), p_{Z_1}(x')) + d^2(C_x, C_{x'})$, and $\alpha$ is isometric as claimed; it is onto by~(\ref{eq_fhi_prop}).

\begin{rems}
\label{rems_exi_bis}%
Let $H$ be a topological group with a continuous action by isometries on a complete \cat0 space $X$. (1)~The above arguments show that the (possibly empty) union of \emph{all} minimal non-empty closed convex $H$\ti invariant subspaces $C$ splits as a product $C\times T$; we call $T$ the \emph{space of components}. (2)~If the action is non-evanescent, then there is a \textbf{canonical} minimal non-empty closed convex $H$\ti invariant subspace $C_0\se X$. Indeed, $T$ is non-empty by Remark~\ref{rem_existence} and bounded by non-evanescence; hence it has a circumcentre $t$ and we let $C_0= C\times\{t\}$.
\end{rems}

\subsection{}
It remains to check that $\alpha$ intertwines the $G$\ti action on $Z$ with the product action on $Z_1\times Z_2$, and actually for the $G_1$ factor this immediately follows from the $G_1$\ti equivariance of $p_{Z_1}$. However, \emph{a priori}, the bijection $\alpha$ transports the $G$\ti action on $Z$ to a $G$\ti action on $Z_1\times Z_2$ of the form
$$g_1 g_2 (p_{Z_1}(x),C_x) = (g_1(g_2\star p_{Z_1}(x)), g_2 C_x)\kern1cm (g_i\in G_i, x\in Z),$$
where we only know that the assignment $G_2\times Z_1\to Z_1$ defined by
$$(g_2, z)\longmapsto g_2\star z \defq p_{Z_1}(g_2 z)$$
determines a well-defined $G_2$\ti action on $Z_1$, which moreover commutes with the $G_1$\ti action. We need to show that $g_2\star = \id_{Z_1}$ for all $g_2\in G_2$. As in the proof of Proposition~\ref{prop_coherent} (for $\teta$), we deduce from Lemma~\ref{lemma_const} that $g_2\star$ is a Clifford translation of $Z_1$; so if it were non-trivial it would preserve a line $\sigma$ in $Z_1$ and in particular fix a point $\sigma(\infty)\in\p X$. Moreover, we would have $g_2^n\star\sigma(0) = \sigma(n\lambda)$ for some $\lambda\neq 0$ and all $n\in {\bf Z}$. Thus for all $g_1\in G_1$ and all $n\in {\bf Z}$ we have
$$d(g_1\sigma(n\lambda), \sigma(n\lambda)) =  d(g_1 g_2^n\star\sigma(0), g_2^n\star\sigma(0)) = d(g_1\sigma(0), \sigma(0))$$
which implies that $g_1$ fixes $\sigma(\infty)$. On the other hand, $G_2$ (not just $G_2\star$) fixes $\sigma(\infty)$ since $\sigma$ lies in $Z_1$. Now there is a $G$\ti fixed point in $\p X$, contradicting the assumption according to Proposition~\ref{prop_eva_proper}.

Being isometric to a product of geodesic spaces, $Z$ is itself geodesic and hence convex in $X$. Likewise, it is closed in $X$ because it is a product of complete spaces (indeed, the uniform structure on $Z_2$ is complete because it is in fact a product uniform structure; alternatively, apply the hands-on argument in~\cite[p.~240]{Bridson-Haefliger}).

\begin{rem}
\label{rem_add}%
It does not follow \emph{a priori} from the above proof that $Z_2$ is minimal. However, we may preface the whole proof by replacing $X$ with the canonical component provided by Remarks~\ref{rems_exi_bis}. Therefore, in the non-evanescent case of Theorem~\ref{thm_splitting}, we obtain in addition that each $Z_i$ is minimal.
\end{rem}

This concludes the proof of Theorem~\ref{thm_splitting}.\hfill\qedsymbol

\section{Induction and its Properties}
\label{sec_suspensions}%
\subsection{}
We begin by defining general ``Pythagorean integrals'' of metric spaces; this is not the only natural integral of metric spaces, see Remark~\ref{rem_join_integral}.

\begin{defi}
\label{defi_int_spaces}%
Let $(\F,\mu)$ be a standard Borel space with a probability measure $\mu$ and let $X$ be a metric space. We denote by $L^2(\F,X)$ the space of all measurable maps (up to null-sets) $f:\F\to X$ with separable range and such that for some (and hence any) $x\in X$ the function $g\mapsto d(f(g), x)$ is in $L^2(\F)$. We endow $L^2(\F,X)$ with the metric defined by
$$d(f, f') \defq \left(\int_\F d^2(f(g), f'(g))\d \mu(g)\right)^{1/2}.$$
\end{defi}

\begin{rems}
(i)~The $L^2$ condition is independent of $x$ by the triangle inequality in $X$ since $\mu$ is finite. (ii)~The triangle inequality in $L^2(\F,X)$ follows by combining the Cauchy-Schwarz inequality with the triangle inequality in $X$.
\end{rems}

Such spaces were considered \emph{e.g.} in~\cite{Korevaar-Schoen}. The following is straightforward:

\begin{lemma}
\label{lemma_ind_top}%
Suppose that $X$ is complete, respectively separable. Then so is $L^2(\F,X)$.\hfill\qedsymbol
\end{lemma}

We now describe geodesic segments in $L^2(\F,X)$; compare~\cite[I.5.3]{Bridson-Haefliger}.

\begin{prop}
\label{prop_geod_proj}%
Let $X$ be a complete metric space, $L^2(\F,X)$ as in Definition~\ref{defi_int_spaces} and $I\se {\bf R}$ any interval. A continuous map $\sigma: I\to L^2(\F,X)$ is a geodesic if and only if there is a mesurable map $\alpha:\F\to{\bf R}_+$ and a collection $\{\sigma^g\}_{g\in\F}$ of geodesics $\sigma^g: \alpha(g)I \to X$ such that
$$\int_\F\alpha(g)^2\d \mu(g) = 1,\kern1cm \sigma(t)(g) = \sigma^g(\alpha(g)t)$$
for all $t\in I$ and $\mu$-a.e. $g\in \F$ (\emph{i.e.} $\alpha$ is a semi-density).
\end{prop}

\begin{proof}
The condition is sufficient; conversely, suppose that $\sigma$ is a geodesic. It suffices to show that there is a dense subset $J\se I$ such that for a.e. $g\in \F$, the map $\sigma(\cdot)(g)$ coincides on $J$ with a linearly reparametrized geodesic $I\to X$. Indeed, if $\alpha(g)$ denotes the reparametrization factor and $\sigma^g$ the corresponding geodesic, we obtain (sufficiency) a geodesic $I\to L^2(\F,X)$ which coincides with $\sigma$ on $J$, hence equals $\sigma$ by continuity. In order to exhibit $J$, it is enough to prove that for every $s<u$ in $I$ and $t=(s+u)/2$ we have
\beq
\label{eq_geod}%
d(\sigma(s)(g), \sigma(u)(g)) = 2d(\sigma(s)(g), \sigma(t)(g)) = 2d(\sigma(t)(g), \sigma(u)(g))
\eeq
for $\mu$\ti a.e $g$. The triangle inequality and $(a+b)^2\leq 2a^2 + 2b^2$ give
$$d^2(\sigma(s)(g), \sigma(u)(g)) \leq 2d^2(\sigma(s)(g), \sigma(t)(g)) + 2d^2(\sigma(t)(g), \sigma(u)(g))$$
$\mu$\ti a.e., with equality if and only if~(\ref{eq_geod}) holds. Integrating, we find
$$d^2(\sigma(s), \sigma(u)) \leq 2d^2(\sigma(s), \sigma(t)) + 2d^2(\sigma(t), \sigma(u))$$
with equality if and only if~(\ref{eq_geod}) holds $\mu$\ti a.e. But equality does hold here since $\sigma$ is geodesic.
\end{proof}

\begin{rem}
The above definitions and the general facts of this section hold more generally for the integral of a measurable field of metric spaces $X_g$ over $\F$. The only difference is that one must fix a choice of a section of base-points in order to define the $L^2$ condition. With this addition, one may also consider infinite measure spaces $\F$. None of this will be needed or used, therefore we leave details to the reader.
\end{rem}

\subsection{}

We now specialize to the \cat0 setting.

\begin{lemma}
\label{lemma_ind_gen}%
Suppose that $X$ is a complete \cat0 space. Then $L^2(\F,X)$ is also a complete \cat0 space.
\end{lemma}

\begin{proof}
The space $L^2(\F,X)$ is geodesic by (the trivial part of) Proposition~\ref{prop_geod_proj}. We need to check inequality~(\ref{eq_CN}) for $x,c,c'$ in $L^2(\F,X)$ and $m$ the midpoint of any geodesic line from $c$ to $c'$. By Proposition~\ref{prop_geod_proj}, $m(g)$ is the midpoint of a geodesic from $c(g)$ to $c'(g)$ for a.e. $g\in\F$. Therefore, the inequality holds pointwise and thus we can integrate it since it is linear in the squares of the distances.
\end{proof}

\begin{exo}
Let $M$ be a Riemannian manifold of finite volume, and denote by $\omega$ the associated volume form. For any $x\in M$, the space of all positive definite symmetric bilinear forms on the tangent space $T_x M$ which also induce $\omega$ is a \cat0 space, since it is isomorphic to the symmetric space $X$ associated to $\mathbf{SL}_n({\bf R})$, where $n=\mathrm{dim}(M)$. Thus, if $(\F, \mu)$ denotes the (normalized) probability space underlying $(M,\omega)$, then the space of ``$L^2$ Riemannian metrics'' on $M$ inducing $\omega$ is isomorphic to $L^2(\F, X)$. Observe that it is endowed with a natural isometric action of the space of volume-preserving diffeomorphisms of $(M,\omega)$ when $M$ is compact.
\end{exo}

\begin{rem}
\label{rem_join_integral}%
It follows from Proposition~\ref{prop_geod_proj} that for $X$ complete \cat0 the boundary of $L^2(\F,X)$ can be easily described as what we call a \emph{join integral}
$$\p L^2(\F,X) \cong \int_\F^* \p X$$
where the right hand side stands for the set of pairs $(\fhi, \alpha)$ consisting of a measurable map $\fhi:\F\to \p X$ and a semi-density $\alpha$ on $\F$; we identify $(\fhi, \alpha)$ with $(\fhi', \alpha')$ when $\alpha = \alpha'$ $\mu$\ti a.e. and $\fhi = \fhi'$ $(\alpha^2\mu)$\ti a.e. We point out to the interested reader that this join integral behaves well with respect to \emph{Tits geometry}; indeed, our definition of the right hand side $\int_\F^*$ makes sense for more general (fields of) spaces replacing $\p X$ and can be endowed with a natural metric by means of integrating the spherical cosine law. We shall not use any of this since the boundary of non-proper spaces contains too little information for our purposes.
\end{rem}

The next proposition establishes that for a \cat0 space $X$, Euclidean rectangles (\emph{cf.} Section~\ref{sec_sandwich}) in $L^2(\F,X)$ decompose as a ``field of parallelograms'' in $X$ over $\F$ (of course, the latter need not be rectangles, as shows even the simplest possible example of the decomposition ${\bf R}^4 = {\bf R}^2\times {\bf R}^2$, which corresponds to an atomic $\mu$).

\begin{prop}
\label{prop_paral_proj}%
Suppose that $X$ is complete \cat0. Let $I\se {\bf R}$ be an interval and $\sigma_1, \sigma_2 : I\to L^2(\F,X)$ two geodesics determining a Euclidean rectangle. Let $\sigma_i^g$ and $\alpha_i$ be as in Proposition~\ref{prop_geod_proj} for $\sigma_i$, $i=1,2$. Then for almost every $g\in\F$ the function $d(\sigma_1^g(\alpha_1(g)t), \sigma_2^g(\alpha_2(g)t))$ is constant on $t\in I$. Moreover, $\alpha_1 = \alpha_2$ $\mu$\ti a.e.
\end{prop}

\begin{proof}
The function $d(\sigma_1(t), \sigma_2(t))$ is constant on $t\in I$. On the other hand, the functions $d(\sigma_1(\alpha_1(g)t), \sigma_2(\alpha_2(g)t))$ are convex and non-negative. Thus the first part of the proposition follows from the general fact that an integral of the squares of a family of convex non-negative functions on $I$ parametrized by a finite measure space is constant if and only if almost every function in the family is constant.

Since $\sigma_1, \sigma_2$ bound a Euclidean rectangle in $L^2(\F,X)$, we may for any two $t_1, t_2\in I$ apply this first part of the proposition to the two geodesics $[\sigma_1(t_1), \sigma_2(t_1)]$, $[\sigma_1(t_2), \sigma_2(t_2)]$. The conclusion is precisely that $\alpha_1(g)=\alpha_2(g)$ holds for a.e. $g\in \F$.
\end{proof}

\subsection{}
\label{sec_bary}%
For any $f\in L^2(\F,X)$, where $X$ is complete \cat0, there is a unique point $x\in X$ minimising $\int_\F d^2(f(g), x)\d \mu(g)$; this point is called the \emph{barycentre} of $f$. Indeed, we may embed $X$ isometrically into $L^2(\F,X)$ by $x\mapsto \psi_x$, $\psi_x(g)\defq x$. Since the image is a closed convex subspace, the barycentre can be defined by the nearest point projection of $f$ to that image. In particular, it follows from this definition that for the barycentres $x,x'$ of $f,f'$ we have
$$d(x,x')\ \leq\ d(f,f').$$
The special case $f'=\psi_{y}$ yields
$$d^2(x,z)\ \leq\ \int_\F d^2(f(g), z)\d \mu(g)\kern1cm(\forall\,z\in X).$$
Actually, the first inequality can be strengthened to
\beq
\label{eq_bary_1}%
d(x,x')\ \leq \ \int_\F d(f(g), f'(g)) \d \mu(g)
\eeq
(though we will only use this in Appendix~A and with $(\F,\mu)$ replaced by a finite set). Indeed, for $0<\epsi<1$ define $x_\epsi, x'_\epsi$ as in Section~\ref{sec_cat}. Integrating~(\ref{eq_Resh}) with $y=f(g)$, $y'=f'(g)$ and using
$$\int_\F\Big( d^2(f(g), x) + d^2(f'(g), x') \Big)\d\mu(g)\ \leq\ \int_\F\Big( d^2(f(g), x_\epsi) + d^2(f'(g), x'_\epsi) \Big)\d\mu(g)$$
yields~(\ref{eq_bary_1}) when $\epsi$ goes to zero. Notice that~(\ref{eq_bary_1}) also reads
\beq
\label{eq_bary}%
d^2(x,x')\ \leq\ d^2(f,f') - \int_\F\Big(d(f(g), f'(g)) - d(x,x')\Big)^2\d \mu(g).
\eeq
More refined inequalities can be found in~\cite[2.5.2]{Korevaar-Schoen}.

\subsection{}
\label{sec_ind_G}%
A particular case of Definition~\ref{defi_int_spaces} arises as follows. Let $G$ be a locally compact second countable group, $\Gamma<G$ a uniform lattice and $X$ a metric space with a $\Gamma$\ti action by isometries. Since $\Gamma$ is cocompact and $G$ second countable, one can find a relatively compact Borel right fundamental domain $\F\se G$ with the property that for any compact $C\se G$ the set
\beq
\label{eq_nice_fund}%
\big\{\eta\in\Gamma : \F\eta\cap C\neq \varnothing \big\}
\eeq
is finite, see \emph{Exercice~12} of~\cite[VII \S2]{BourbakiINT78} (just mind that the scholar of Nancago considers \emph{left} fundamental domains). \emph{We shall from now on agree to consider only such domains; we further endow $G$ with a Haar measure $\mu$ normalized by $\mu(\F)=1$ and write $\d g$ for $\d\mu(g)$.}

\begin{rem}
We assumed $G$ second countable so that its Borel structure is standard; in addition, we will often assume $X$ separable. We do however prove the theorems of the introduction in the full generality stated there by showing in due time how to reduce to the current assumptions.
\end{rem}

\begin{defi}
\label{defi_induction}%
We denote by $L^{[2]}(G,X)^\Gamma$ the space of all measurable maps (up to null-sets) $f:G\to X$ with separable range and such that (i)~for all $g\in G$, $\gamma\in \Gamma$ one has $f(g\gamma^{-1})=\gamma f(g)$ and (ii)~for some (hence any) $x\in X$, the function $g\mapsto d(f(g), x)$ is in $L^2(\F)$. We endow $L^{[2]}(G,X)^\Gamma$ with the metric defined by
\beq
\label{eq_def_dist}%
d(f, f') \defq \left(\int_{G/\Gamma} d^2(f(g), f'(g))\d g\right)^{1/2}.
\eeq
\end{defi}

There is a canonical isometry
\beq
\label{eq_2ind}%
L^{[2]}(G,X)^\Gamma\cong L^2(\F,X)
\eeq
given by restriction to $\F$. Moreover, the choice of $\F$ is equivalent to the choice of a Borel map $\chi:G\to \Gamma$ such that
\beq
\label{eq_dot_action}%
\chi^{-1}(e)=\F, \kern1cm\chi(g\gamma^{-1})=\gamma\chi(g)\kern1cm(\forall\,\gamma\in\Gamma, \mbox{a.e.}\ g\in G).
\eeq
The isomorphism $\F\to G/\Gamma$ induces a $G$\ti action on $\F$ which is described by the rule $h.g=h g\chi(h g)$ (the dot notation emphasizes the difference between the two actions). Moreover, this action is measure-preserving since the existence of $\Gamma$ forces $G$ to be unimodular. Observe that the inverse to the restriction map in~(\ref{eq_2ind}) consists in extending $f\in L^2(\F,X)$ to a map
\beq
\label{eq_extend}%
f_\mathrm{ext}:G\to X,\kern1cm f_\mathrm{ext}(g) \defq \chi(g)f(g\chi(g)).
\eeq
We will abuse notation in omitting the subscript `ext'.

\begin{lemma}
\label{lemma_ind_cont}%
Assume that $X$ is complete and separable. Then there is a well-defined continuous $G$\ti action by isometries on $L^{[2]}(G,X)^\Gamma$ defined by $(hf)(g) \defq f(h^{-1}g)$, where $g,h\in G$.
\end{lemma}

In other words,~(\ref{eq_extend}) and~(\ref{eq_dot_action}) show that the corresponding $G$\ti action on $L^2(\F,X)$ is
\beq
\label{eq_act}%
(h.f)(g) = \chi(h^{-1}g)f(h^{-1}.g),
\eeq
and that the latter is well-defined. We call this the \emph{induced $G$\ti action}.

\begin{proof}[Proof of the lemma]
First we need to show that $hf$ is still in $L^{[2]}(G,X)^\Gamma$. Since $h^{-1} \F$ is relatively compact, the finiteness of~(\ref{eq_nice_fund}) guarantees that there are $\gamma_1, \ldots,\gamma_k$ in $\Gamma$ such that the union of the $\F\gamma_i$ covers $h^{-1}\F$. Fix a base-point $x\in X$. Then
$$\int_\F d^2(hf(g), x)\d g\ =\ \int_{h^{-1}\F} d^2(f(g), x)\d g\ \leq\ \sum_{i=1}^k \int_{\F \gamma_i} d^2(f(g), x)\d g.$$
But in view of $f(g\gamma_i) = \gamma_i^{-1}f(g)$, each term
$$\int_{\F \gamma_i} d^2(f(g), x)\d g\ =\ \int_{\F} d^2(\gamma_i^{-1}f(g), x)\d g\ =\ \int_{\F} d^2(f(g), \gamma_i x)\d g$$
is finite since $f$ is in $L^{[2]}(G,X)^\Gamma$. This action preserves the distance~(\ref{eq_def_dist}). \emph{Ad} continuity: By Lemma~\ref{lemma_cont_Borel}, it is enough to show that the map $h\mapsto h f$ is measurable for all $f$ in $L^{[2]}(G,X)^\Gamma$. This follows from the fact that the map $G\times G\to X$, $(g,h)\mapsto f(h^{-1}g)$ is measurable. (Alternatively, for $X$ \cat0, one can also show that continuous functions are dense using barycentres weighted by continuous approximate units on $G$.)
\end{proof}

\subsection{}
\label{sec_ind_eva}%
Some properties of the $\Gamma$\ti action on $X$ are trivially equivalent to the corresponding property for the $G$\ti action on $L^{[2]}(G,X)^\Gamma$; for instance, the existence of fixed points. Evanescence is more subtle; a very simple instance is when $X$ is a Hilbert space, in which case one has the following standard theorem: if a unitary representation of the cocompact lattice $\Gamma$ does not weakly contain the trivial representation, then the induced $G$\ti representation does not either (this follows \emph{e.g.} from the topological Frobenius reciprocity of~\cite{Blanc}). The following is a geometric generalization:

\begin{thm}
\label{thm_ind_eva}%
Let $G$ be a locally compact second countable group, $\Gamma<G$ a uniform lattice and $X$ a complete separable \cat0 space with a $\Gamma$\ti action by isometries. If the $G$\ti action on $L^{[2]}(G,X)^\Gamma$ is evanescent, then the $\Gamma$\ti action on $X$ is evanescent.
\end{thm}

Observe that one cannot reduce this to a statement about asymptotic cones through Remark~\ref{rem_ultra}; indeed, even when $X=\mathbf{R}$, an asymptotic cone on $L^{[2]}(G,X)^\Gamma$ is already a so-called \emph{non-standard hull} of the Hilbert space $L^{[2]}(G,X)^\Gamma$~-- whilst $X$ is its own asymptotic cone.

Instead, the general idea is to project an evanescent sequence of functions to the space constant functions. This does not quite work since one needs to spread out the domain of these functions beyond any fundamental domain in order to capture generators of $\Gamma$. The following proof is particularly simple thanks to Lemma~\ref{lemma_ind_eva}; another argument is given in Appendix~B for certain non-uniform lattices.

\begin{proof}[Proof of Theorem~\ref{thm_ind_eva}]
The finiteness of the set in~(\ref{eq_nice_fund}) implies in particular that for every relatively compact Borel set $\E\se G$ of positive measure and any $f$ in $L^{[2]}(G,X)^\Gamma$ the restriction $f|_\E$ is in $L^2(\E,X)$ (as in Lemma~\ref{lemma_ind_cont}). Therefore we may define a point $x_{f,\E}\in X$ by taking the barycentre of $f|_\E$; that is, $x_{f,\E}\in X$ minimises
$$\Delta_\E(f) \ \defq\ \inf_{x\in X}\int_\E d^2(x,f(g))\d g.$$
For every $x\in X$, define the element $\psi_x$ of $L^{[2]}(G,X)^\Gamma$ by $\psi_x(g) \defq \chi(g)x$ (extending the notation of Section~\ref{sec_bary}); finally, write $f_\E \defq \psi_{x_{f,\E}}$.

\smallskip

Fix a finite set $S\se \Gamma$. We shall produce an unbounded $S$\ti evanescent set $T_S\se X$ for every unbounded $Q$\ti evanescent set $T$ in $L^{[2]}(G,X)^\Gamma$ with $Q\se G$ a suitable compact set. This shows that weak evanescence of the $G$\ti action on $L^{[2]}(G,X)^\Gamma$ implies weak evanescence for the $\Gamma$\ti action on $X$; the statement for evanescence follows by Proposition~\ref{prop_ceva} and Remark~\ref{rem_ceva}.

\smallskip

Define a relatively compact set $\E$ and choose a compact set $Q$ with
$$\E\ \defq\ \bigcup_{\gamma\in S} \gamma^{-1} \F \gamma\ \cup \F, \kern1cm Q \supseteq \F\E^{-1} \cup S.$$

\begin{lemma}
\label{lemma_ind_eva}%
The function $\Delta_\E$ is bounded on each $Q$\ti evanescent set in $L^{[2]}(G,X)^\Gamma$.
\end{lemma}

\begin{proof}[Proof of the lemma]
Let $T$ be a $Q$\ti evanescent set. There is $K$ such that $d^2(f, hf)\leq K$ for all $h\in \F\E^{-1}$ and all $f\in T$. Let $f\in T$; we have
$$\Delta_\E(f) \ =\ \int_\E d^2(x_{f,\E},f(g))\d g\ \leq\ \int_\F\int_\E d^2(f(\bar g),f(g))\d g \d\bar g$$
by the choice of $x_{f,\E}$. After the change of variable $h=\bar g g^{-1}$ this is
$$\int_\F\int_{\bar g \E^{-1}} d^2(f(\bar g), f(h^{-1}\bar g))\d h \d\bar g\ \leq\ \int_\F\int_{\F\E^{-1}} d^2(f(\bar g), f(h^{-1}\bar g))\d h \d\bar g,$$
which is just
$$\int_{\F\E^{-1}}\int_\F d^2(f(\bar g), hf(\bar g))\d \bar g \d h\ =\ \int_{\F\E^{-1}} d^2(f, hf)\d h.$$
This is bounded in terms of $K$ and the measure of $\F\E^{-1}$.
\end{proof}

Choose now an unbounded $Q$\ti evanescent set $T$ in $L^{[2]}(G,X)^\Gamma$ and consider the subset $T_S = \{x_{f,\E} : f\in T \}$ of $X$. We contend that $T_S$ is $S$\ti evanescent. First, we check that $T_S$ is unbounded: Indeed, since $\F\se\E$, for all $f\in T$
\beq
\label{eq_ind_eva}%
d^2(f, f_\E)\ =\ \int_\F d^2(f(g), x_{f,\E})\d g\ \leq \Delta_\E(f)
\eeq
is bounded by Lemma~\ref{lemma_ind_eva}, and therefore this first claim follows from the estimate
$$d(x_{f, \E}, x_{f', \E}) \ =\ d(f_\E, f'_\E) \ \geq\ d(f,f') - d(f, f_\E) - d(f', f'_\E)$$
for $f, f'\in T$ since $T$ is unbounded. We now need to estimate $d(x_{f,\E}, \gamma x_{f,\E})$ uniformly over $\gamma\in S$ and $f\in T$. Write
$$d(x_{f,\E}, \gamma x_{f,\E}) \ =\ d(f_\E, \psi_{\gamma x_{f,\E}}) \ \leq\ d(f_\E, f) + d(f,\gamma f) + d(\gamma f,\psi_{\gamma x_{f,\E}}).$$
The first term is taken care of by~(\ref{eq_ind_eva}) and Lemma~\ref{lemma_ind_eva}, whilst the second is bounded since $T$ is $Q$\ti evanescent and $S\se Q$. As for the last term, we have
\begin{multline*}
d^2(\gamma f,\psi_{\gamma x_{f,\E}}) \ =\ \int_\F d^2(f(\gamma^{-1}g), \gamma x_{f,\E})\d g \ =\ \int_\F d^2(\gamma^{-1}(f(\gamma^{-1}g)), x_{f,\E})\d g\\
=\ \int_\F d^2(f(\gamma^{-1}g\gamma), x_{f,\E})\d g \ =\ \int_{\gamma^{-1}\F\gamma} d^2(f(g), x_{f,\E})\d g\\
\leq\ \int_\E d^2(f(g), x_{f,\E})\d g\ =\Delta_\E(f),
\end{multline*}
so we are done by Lemma~\ref{lemma_ind_eva}.
\end{proof}

\section{Superrigidity}
\label{sec_proofs}%
Throughout this section, when considering irreducible lattices $\Gamma<G=G_1\times \cdots\times G_n$, we shall always assume $n\geq 2$. This is indeed not a restriction, because for $n= 1$ the definition of irreducibility implies $G=\Gamma$ and there is nothing to prove.

\subsection{}
\label{sec_super_intermediate}%
We now give, in a slightly cumbersome formulation, a key intermediate statement to which various superrigidity statements will be reduced. A map is called \emph{totally geodesic} if it takes geodesic segments to (possibly reparametrized) geodesic segments. Recall that the separability assumptions, made as a matter of convenience, will be disposed of in due time.

\begin{thm}
\label{thm_superrigidity2}%
Let $\Gamma$ be an irreducible uniform lattice in a product $G=G_1\times \cdots\times G_n$ of locally compact second countable groups and let $X$ be a complete separable \cat0 space with a non-evanescent $\Gamma$\ti action by isometries.

Then there is a (canonical and minimal) non-empty closed convex $G$\ti invariant subspace $Z\se L^{[2]}(G,X)^\Gamma$ which splits isometrically and $G$\ti equivariantly as a product $Z_1\times \cdots\times Z_n$ of minimal $G_i$\ti spaces $Z_i$. Moreover, $Z$ consists of continuous functions and the evaluation map $\psi: Z\to X$ at $e\in G$ is a totally geodesic $\Gamma$\ti equivariant Lipschitz map. The restriction of $\psi$ to every copy of each $Z_i$ is isometric. In particular, if for some $i$ the set $\C_i \defq \p(\psi(Z_i))\se \p X$ is non-empty, then $\Gamma$ preserves $\C_i$ and the $\Gamma$\ti action on $\C_i$ extends continuously to a $G$\ti action which factors through $G\to G_i$.
\end{thm}

As usual, $f\in L^{[2]}(G,X)^\Gamma$ is really a function \emph{class} and is said continuous if it contains a continuous representative. We shall use the following criterion, which is readily checked using Fubini-Lebesgue.

\begin{lemma}
\label{lemma_crit_cont}%
Let $f:G\to X$ be a measurable function from a locally compact second countable group $G$ to a separable complete metric space $X$. Assume that for all $\epsi>0$ and all $g\in G$ there is a neighbourhood $U$ of $g$ such that $\smash{\esssup\limits_{h,h'\in U}}\, d(f(h), f(h'))<\epsi$. Then $f$ agrees almost everywhere with a continuous function.\hfill\qedsymbol
\end{lemma}

\begin{proof}[Proof of Theorem~\ref{thm_superrigidity2}]
Recall that $L^{[2]}(G,X)^\Gamma$ is the complete separable \cat0 space with continuous $G$\ti action granted by Definition~\ref{defi_induction} and Lemmata~\ref{lemma_ind_top},~\ref{lemma_ind_gen} and~\ref{lemma_ind_cont}. Theorem~\ref{thm_ind_eva} shows that the $G$\ti action is not evanescent, and thus the splitting Theorem~\ref{thm_splitting} implies that there is a closed convex $G$\ti invariant subspace $Z\se L^{[2]}(G,X)^\Gamma$ which splits isometrically and $G$\ti equivariantly as a product $Z_1\times \cdots\times Z_n$ of $G_i$\ti spaces $Z_i$ endowed with the product action. Recall from Remark~\ref{rem_add} that $Z_i$ can be assumed to be a minimal non-empty closed convex $G_i$\ti space.

We may, and shall from now on, identify each $Z_i$ with a subspace of $Z$. For every $y\in Z$, there is a unique copy of $Z_i$ in $Z$ containing $y$; we denote by $Z_i^y$ this $G_i$\ti invariant closed convex subset $Z_i^y\se Z$. We claim
\beq
\label{eq_psi_isom}%
d(f(g), f'(g))\ =\ d(f,f') \kern1cm \forall\,y\in Z, \forall\, f,f'\in Z_i^y, \text{ a.e. } g\in G.
\eeq
Indeed, fix a Borel fundamental domain $\F\se G$ for $\Gamma$ as in Section~\ref{sec_ind_G}. We shall make use of the identification~(\ref{eq_2ind}) and of the notation~(\ref{eq_act}). It suffices to prove the claim for a.e. $g\in\F$. Let $I=[0,d(f,f')]$ and let $\sigma_1:I\to Z_i^y$ be the geodesic from $f$ to $f'$. If $h$ is any element of $G_j$ for $j\neq i$, the splitting of $Z$ shows that $\sigma_1$ and $\sigma_2 \defq h^{-1}\sigma$ determine a Euclidean rectangle in $Z$. Let $\sigma_i^g$ and $\alpha_i$ be as in Proposition~\ref{prop_geod_proj} for $\sigma_i$, $i=1,2$. Then Proposition~\ref{prop_paral_proj} shows that $\alpha_1 = \alpha_2$ a.e. In other words, since $\alpha_2(g) = \alpha_1(h.g)$, the function (class) $\alpha_1:\F\to{\bf R}_+$ is $G_j$\ti invariant for all $j\neq i$. Since the projection of $\Gamma$ to $G_i$ is dense, the subproduct $\prod_{j\neq i} G_j$ acts ergodically on $G/\Gamma$. Therefore, $\alpha_1$ (and hence also $\alpha_2$) is constant; this constant is one by $\int_\F\alpha_1(g)^2\d g = 1$. The claim now follows since $f(g)$ and $f'(g)$ are the endpoints of $\sigma_1^g$.

\smallskip

Next, we recall (for any $a_i\in{\bf R}$) the inequality $\sum_{i=1}^n a_i\leq \big(n\sum_{i=1}^n a_i^2\big)^{1/2}$; this allows us to bound the distance between any $f,f'\in Z$ by applying~(\ref{eq_psi_isom}) to each factor $Z_i$, obtaining
\beq
\label{eq_Lip}%
d(f(g), f'(g)) \leq \sqrt{n}\, d(f,f')\kern1cm \forall\,f,f'\in Z, \text{ a.e. } g\in G.
\eeq
Since $Z$ is $G$\ti invariant, it follows that for any compact neighbourhood $U$ of $g\in G$
$$\esssup_{h,h'\in U} d(f(h), f(h'))\ \leq\ \sqrt{n}\, \sup_{h,h'\in U} d(f,(hh'^{-1})f)\kern1cm \forall\,f\in Z$$
which goes to zero as $U\to g$ by continuity of the $G$\ti action (Lemma~\ref{lemma_ind_cont}). This shows that every $f\in Z$ is continuous (Lemma~\ref{lemma_crit_cont}).

\medskip

We may now define a map $\psi: Z\to X$ by $\psi(f) = f(e)$. This map is $\Gamma$\ti equivariant by definition. For every $y\in Z$ and each $i$, the restriction of $\psi$ to $Z_i^y\to X$ is isometric because of~(\ref{eq_psi_isom}). Further, $\psi$ is $\sqrt{n}$\ti Lipschitz by~(\ref{eq_Lip}) and totally geodesic by Proposition~\ref{prop_geod_proj}. Assume, for some $i$, that $\p Z_i\neq \varnothing$; set $X_i \defq \psi(Z_i)$. Notice that $X_i$, being isometric to $Z_i$, is complete, hence closed, and convex. In particular, $\C_i \defq \p X_i$ is closed in $\p X$~\cite[p.~266]{Bridson-Haefliger} and the isometry $\psi|_{Z_i}:Z_i\to X_i$ induces a homeomorphism $\p Z_i \cong \C_i$. The $G$\ti action on $Z$ being the product action, $\p Z_i$ is $G$\ti invariant and the (continuous) $G$\ti action on it factors through $G\to G_i$. Summing up, we obtained a $\Gamma$\ti equivariant homeomorphism $\p Z_i\cong \C_i$, finishing the proof.
\end{proof}

\subsection{}
\label{sec_decomposable}%
In this subsection we analyse the notion of indecomposability introduced in Definition~\ref{defi_decomposable} (for proper spaces). The two obvious obstructions to indecomposability of a group $L$ are (i)~a $L$\ti fixed point at infinity (when $L\neq 1$) and (ii)~a $L$\ti invariant splitting $X=X_1\times X_2$ with unbounded factors; indeed in both cases we obtain closed invariant sets at infinity on which $L$ does not act faithfully. It is not clear to us to what extent the topological part of Definition~\ref{defi_decomposable} is really an additional restriction in our setting.

\smallskip

Let $X$ be a proper \cat0 space, $L<\isom(X)$ a subgroup and $\C\se \p X$ a non-empty closed $L$\ti invariant subset; write $J$ for the stabiliser of $\C$ in $\isom(X)$. Recall that the topology of uniform convergence coincides on $M\defq \homeo(\C)$ with the compact-open topology and turns $M$ into a topological group; the natural homomorphism $\iota:J\to M$ is continuous (compare Section~\ref{sec_boundary}). Moreover $M$ is polish (hence Baire) because $\C$ is metrizable. Notice that  $J$ is closed in $\isom(X)$ and hence contains $\bbb L$.

\begin{rem}
\label{rem_indec_closed}%
The conditions of Definition~\ref{defi_decomposable} hold if and only if $\iota|_{\bbb L}$ is a topological isomorphism onto its image. Indeed, indecomposability implies that $\iota|_{\bbb L}$ is a continuous group isomorphism from a locally compact second countable group onto a Baire group and thus is a topological isomorphism by the usual Baire category argument. The converse (which we do not use) holds since any locally compact subgroup of a topological group is closed (Corollaire~2 in~\cite[II \S3 \no3]{BourbakiTG34}).
\end{rem}

In the two general instances where we verify indecomposability (Lemmata~\ref{lemma_Zariski} and~\ref{lem_hyperbolic}), we obtain the \emph{a priori} stronger statement for the coarser topology of pointwise convergence by the following criterion (mind that $M$ is not a topological group for this topology).

\begin{lemma}
\label{lem_crit_indec}%
Let $\bbb L$ be a locally compact second countable group with a faithful continuous action $\bbb L\times \C \to \C$ on a compact Hausdorff topological space $\C$. The associated injective homomorphism $\bbb L\to M\defq \homeo(\C)$ is a topological isomorphism onto its image endowed with the topology of pointwise convergence if and only if for any sequence $\{h_n\}$ tending to infinity in $\bbb L$ there is $\xi\in \C$ with  $h_n \xi \nrightarrow \xi$.
\end{lemma}

\begin{proof}
Let $\{g_\alpha\}$ be any net of $\bbb L$ converging to $g\in\bbb L$ for the pointwise topology on $M$ but not in the $\bbb L$\ti topology. Since $\bbb L$ is locally compact second countable and its action on $\C$ continuous, we obtain a sequence $\{\ell_n\}$ converging to infinity in $\bbb L$ such that $\ell_n\xi\to g\xi$ for all $\xi\in C$. The criterion applied to $h_n\defq g^{-1} \ell_n$ yields a contradiction. The converse follows from the continuity of the action.
\end{proof}

We now verify that indecomposability generalizes indeed Zariski-density for subgroups of adjoint simple algebraic groups.

\begin{lemma}
\label{lemma_Zariski}%
Let $k$ be a local field and $\mathbf{H}$ a connected adjoint $k$\ti simple $k$\ti group. Let $X$ be the symmetric space, respectively the Bruhat-Tits building, associated to $G=\mathbf{H}(k)$ according to whether $k$ is Archimedean or not. Then any Zariski-dense subgroup $L<G$ is indecomposable.
\end{lemma}

The following proof was kindly provided by J.-F.~Quint.

\begin{proof}
Let $\C\se\p X$ be a non-empty $L$\ti invariant closed subset and $H_\C$ its pointwise stabiliser in $G$. Since $H_\C$ is an intersection of parabolic subgroups, it is an algebraic subgroup; but $L$ normalizes it because it preserves $\C$, so by Zariski-density $H_\C$ is normal in $G$ and hence trivial.

Let now $\{h_n\}$ be a sequence tending to infinity in $\bbb L$. Since parabolic subgroups are cocompact, there is no loss of generality in supposing that $\C$ is contained in a single $G$\ti orbit; thus $\C$ can be considered as a Zariski-dense subset of $G/Q$ for some parabolic $Q\neq G$. We may fix an irreducible representation $r$ of $G$ on a $k$\ti vector space $V$ (of dimension, say, $m$) such that $Q$ is the stabiliser of some line $U\se V$, see~\cite{Tits71}. Suppose for now that $k={\bf R}$; the Cartan decomposition gives $r(h_n)= c'_n a_n c_n$ for some $a_n$ in a positive Weyl chamber and $c_n, c'_n\in\mathbf{SO}(m)$. Upon passing to a subsequence, we may assume that $c_n, c'_n$ converge to some $c, c'$. Since $r(h_n)$ tends to infinity in $\mathbf{SL}_m(k)$, the sequence $a_n/\|a_n\|$ (where $\|a_n\|$ is the operator norm) converges, after possibly extracting a subsequence, to a linear map onto a proper subspace $0\neq W \subsetneqq V$, with kernel $0\neq Y\subsetneqq V$. Since $r$ is irreducible, $G U \nsubseteq c' W\cup c^{-1} Y$ so that by Zariski-density there is $\xi\in\C$ not in $c' W\cup c^{-1}Y$. However, for any line $U'\nsubseteq c^{-1}Y$, $r(h_n)U'$ tends to a line in $c' W$. Thus $h_n\xi\nrightarrow \xi$. For general $k$, the same argument applies to the corresponding Cartan decomposition of $\mathbf{GL}_m$ over $k$.
\end{proof}

For the sake of illustration, we now take a brief look at the context of negative curvature (even though the methods proposed in~\cite{Monod-Shalom1},\cite{Monod-ShalomCRAS} are more powerful in that case). Recall that if the proper \cat0 space $X$ is Gromov-hyperbolic, \emph{e.g.} \cat{-1}, an action on $X$ is said \emph{elementary} if it fixes a point in $\bbb X$ or stabilises a pair of points in $\p X$.

\begin{lemma}
\label{lem_hyperbolic}%
Let $\Gamma$ be a group with a non-elementary action by isometries on a Gromov-hyperbolic proper \cat0 space $X$. Then, upon possibly replacing $X$ by a non-empty closed convex $\Gamma$\ti invariant subspace, the $\Gamma$\ti action on $X$ is indecomposable.
\end{lemma}

(The $\Gamma$\ti action is called indecomposable if the image of $\Gamma$ in $\isom(X)$ is an indecomposable subgroup.)

\begin{proof}
By Remark~\ref{rem_existence}, there is a minimal non-empty closed convex $\Gamma$\ti invariant subspace $Y\se X$; observe that the $\Gamma$\ti action on $Y$ is still non-elementary. Thus we may assume $X$ minimal. Recall that there is a continuous $\isom(X)$\ti equivariant map $C$ from distinct triples in $\p X$ to $X$ (indeed, for any distinct $\xi_i\in \p X$, the sum $\beta = \beta_{\xi_1} + \beta_{\xi_2} + \beta_{\xi_3}$ of Busemann functions is bounded below and proper as shown by a comparison argument \emph{e.g.} using the \emph{arbre approximatif} of Th\'eor\`eme~12~(ii) in~\cite[2 \S2]{Ghys-Harpe}; therefore, one can take for $C(\xi_i)$ the circumcentre of its min-set). Let now $\C\se\p X$ be a non-empty closed $\Gamma$\ti invariant set. By non-elementarity, $\C$ contains at least three distinct points $\xi_i$. The pointwise stabiliser $K$ of $\C$ in $\isom(X)$ is therefore compact since it fixes $C(\xi_i)$; thus, the closed convex subspace $X^K$ is non-empty. Since $\Gamma$ normalizes $K$, it preserves $X^K$; hence, $X^K = X$ and thus $K$ is trivial.

Let now $\{h_n\}$ be any sequence in the stabiliser $H<\isom(X)$ of $\C$ such that $h_n\xi \to\xi$ for all $\xi\in \C$. It remains only to show that the sequence $\{h_n\}$ is bounded in $H$, or equivalently in $\isom(X)$. This follows since $h_n C(\xi_i)$ tends to $C(\xi_i)$.
\end{proof}

\subsection{}
We now proceed to prove Theorem~\ref{thm_boundary} and Corollary~\ref{cor_superrigidity}.

\begin{prop}
\label{prop_s_red}%
It is enough for Theorem~\ref{thm_boundary} to consider the case where $G$ is second countable.
\end{prop}

\begin{proof}
Based on a classical argument~\cite{Kakutani-Kodaira}, one has:

\begin{lemma}
\label{lemma_KK}%
Let $J$ be a locally compact $\sigma$\ti compact group and $V\se J$ a neighbourhood of $e\in J$. Then there exists a compact normal subgroup $K\lhd J$ contained in $V$ such that $J/K$ is second countable.
\end{lemma}

\begin{proof}[Proof of the lemma]
Let $\fhi$ be a non-negative continuous function on $J$ that is supported on $V$ and such that $\fhi(e)=1$. Then the smallest closed $J$\ti invariant subspace $\M$ of $L^2(J)$ containing $\fhi$ is separable, see Satz~5 in~\cite{Kakutani-Kodaira}. If we let $K$ be the kernel of the $J$\ti representation on $\M$, it follows as in~\cite{Kakutani-Kodaira} that $J/K$ is second countable because, by $\sigma$\ti compactness and an application of Baire's theorem, $J/K$ is topologically isomorphic to its image in the unitary group of $\M$ endowed with the strong operator topology. The choice of $\fhi$ guarantees that $g\fhi\neq \fhi$ for all $g\notin V$, so that $K\se V$.
\end{proof}

Since $\Gamma$ is discrete in $G$, there are identity neighbourhoods $V_i\se G_i$ such that $\Gamma$ meets $V_1\times\cdots\times V_n$ trivially. Let $K_i\lhd G_i$ be as in the lemma with $K_i\se V_i$ and set $K=K_1\times \cdots \times K_n$. Then $G^*=G/K$ is second countable and the canonical image $\Gamma^*$ of $\Gamma$ in $G^*$ is still an irreducible uniform lattice with respect to the product structure given by the factors $G_i/K_i$. But the choice of $V_i$ ensures that the natural map $\Gamma\to \Gamma^*$ is an isomorphism; thus, if we have Theorem~\ref{thm_boundary} for the second countable case, we can apply it to $G^*$ and the general case follows.
\end{proof}

\begin{proof}[End of proof of Theorem~\ref{thm_boundary}]
Suppose first that the $\Gamma$\ti action on $X$ is evanescent. Since $X$ is proper, it follows from Proposition~\ref{prop_eva_proper} that there is a $\Gamma$\ti fixed point $\xi\in\p X$. In that case the conclusion of the theorem holds trivially for $\C=\{\xi\}$.

Thus we may suppose that the action is non-evanescent; in addition, we may assume $G$ second countable (Proposition~\ref{prop_s_red}) and $X$ is separable since it is a proper metric space. Therefore Theorem~\ref{thm_superrigidity2} applies; we shall use its notation. If every $Z_i$ is bounded, then $Z$ is reduced to a point by minimality. The image of this point under $\psi$ is a $\Gamma$\ti fixed point in $X$, contradicting unboundedness. Hence we may assume that there is an index $i$ such that $Z_i$ is unbounded. Since $Z_i$ is isometric to the unbounded closed convex subspace $\psi(Z_i)$ of the \emph{proper} space $X$, we deduce that its boundary is non-empty. Now the conclusion of Theorem~\ref{thm_superrigidity2} is exactly what was to be shown.
\end{proof}

\begin{proof}[Proof of Corollary~\ref{cor_superrigidity}]
Keep the notation of the corollary, set $L\defq \tau(\Gamma)$ and let $\C$ be as in Theorem~\ref{thm_boundary}. The conclusion of that theorem implies that the composed map
$$\Gamma \lra \bbb L \xrightarrow{\ \iota\ } M\defq \homeo(\C),$$
with $M$ topologized as in Section~\ref{sec_decomposable}, extends to a continuous homomorphism $\hat\tau:G\to M$ factoring through one of the canonical projections $G\to G_i$. Since the image of $\Gamma$ in $G_i$ is dense and $\iota({\bbb L})$ closed, $\hat\tau(G)\se \iota(\bbb L)$. By Remark~\ref{rem_indec_closed}, $\hat\tau$ lifts to a continuous homomorphism $\wt{\tau}:G\to \bbb L$ factoring through $G_i$ and extending $\tau$. This proves the theorem since $\bbb L<H$.
\end{proof}

\begin{proof}[Proof of Margulis' Theorem~\ref{thm_Margulis} for $\Gamma$ cocompact]
By Lemma~\ref{lemma_Zariski}, we may apply Corollary~\ref{cor_superrigidity}.
\end{proof}

\subsection{}
\label{sec_non-proper}%
With Theorem~\ref{thm_non-proper} in view, we now turn to complete \cat0 spaces $X$ that are not assumed proper and analyse \emph{reduced} subgroups of $\isom(X)$ following Definition~\ref{defi_reduced}. In the beginning of Section~\ref{sec_decomposable} we mentioned two immediate restrictions following from indecomposability; the analogous two restrictions hold in the present setting aswell, as shown by~(i) and~(ii) in the following lemma.

\begin{lemma}
\label{lemma_reduced}%
Let $X$ be a complete \cat0 space and $L<\isom(X)$ an unbounded reduced subgroup. Then:

\noindent
(i)~The action is not evanescent.

\noindent
(ii)~There is no $L$\ti invariant splitting $X=X_1\times X_2$ with non-trivial factors.

\noindent
(iii)~If $L$ is countable, then $X$ is separable.

\noindent
(iv)~$X$ has no non-trivial Clifford translation unless it is isometric to a Hilbert space.
\end{lemma}

\begin{proof}
(ii)~Assume that the splitting $X=X_1\times X_2$ is preserved by $L$. Since $L$ is unbounded, one of the factors, say $X_2$, is unbounded. For any $x\in X_1$ the set $Y=\{x\}\times X_2$ enjoys the property of Definition~\ref{defi_reduced}. Therefore $Y=X$ and hence $X_1$ is trivial.

(iv)~This follows from~(ii) in view of the canonical splitting of $X$ recalled at the end of Section~\ref{sec_sandwich}.

(i)~Suppose for a contradiction that $T\se X$ is an unbounded evanescent set. Since for every $\gamma\in \Gamma$ the displacement function $x\mapsto d(\gamma x, x)$ is continuous and convex, there is no loss of generality in assuming $T$ closed and convex upon replacing it by $[T]$. Then Definition~\ref{defi_reduced} implies $T=X$. It follows now from this same condition that the unbounded \cat0 space $X$ contains no unbounded closed convex subset $Y\subsetneqq X$ at all. In particular, $X$ cannot be a Hilbert space and thus, by~(iv), $X$ does not admit any non-trivial Clifford translation. Since $\Gamma$ acts non-trivially, there is $\gamma\in\Gamma$ whose displacement function is not identically zero; however, since $X=T$, the supremum $C=\sup_{x\in X} d(\gamma x, x)$ is finite. If the displacement length were constant, $\gamma$ would be a non-trivial Clifford translation. Thus for some $0<C'<C$ the closed convex set
$$Y \defq \big\{ x\in X : d(\gamma x, x) \leq C' \big\} \ \subsetneqq X$$
is non-empty, hence bounded. Since $\gamma Y = Y$, the circumcentre $y\in Y$ is $\gamma$\ti fixed. Let $\ro$ be the circumradius of $Y$. Since we assumed $T$ unbounded, there is a sequence $\{x_n\}$ in $X$ with $d(x_n, y)\to\infty$. Let $y_n$ be the point of $[y, x_n]$ at distance $\ro+1$ of $y$; then a comparison argument shows that $d(\gamma y_n, y_n)$ tends to zero since $d(\gamma x_n, x_n)\leq C$. Thus, for $n$ large enough, $y_n$ is in $Y$, a contradiction proving~(i).

(iii)~Let $Y$ be the closed convex hull of some $L$\ti orbit. Since $Y$ is unbounded and $L$\ti invariant, $X=Y$. Thus it remains to show that $Y$ is separable. This is indeed the case: If $A_0\se X$ is any countable set (such as a $L$\ti orbit) and $A_{n+1}$ denotes the set obtained by adjoining to $A_n$ all midpoints of pairs in $A_n$, then the union $\bigcup_{n\geq 0}A_n$ contains the midpoint of any of its pairs, and thus is dense in $[A_0]$.
\end{proof}

We can now complete the main case of Theorem~\ref{thm_non-proper}.

\begin{proof}[End of proof of Theorem~\ref{thm_non-proper}, non-Hilbertian case]
Observe first that $\Gamma$ is countable since it is a lattice in a locally compact $\sigma$\ti compact group. The argument of Proposition~\ref{prop_s_red} shows that we may assume $G$ second countable. Moreover, since $\tau(\Gamma)$ is unbounded and reduced, Lemma~\ref{lemma_reduced} implies that its action is non-evanescent and that $X$ is separable.

Thus, we may apply Theorem~\ref{thm_superrigidity2}; we keep its notation and denote by $z_0\in Z$ the point common to all $Z_i$. We know, as in the proof of Theorem~\ref{thm_boundary}, that at least some $Z_i$ is unbounded; say $Z_n$. Recall that $\psi(Z_n)$ is a closed convex subset of $X$ since it is the isometric image of a complete convex set. Fix $\gamma\in\Gamma$ and write $\gamma = \gamma_n \gamma'_n$ for the decomposition along $G= G_n\times G'_n$, where $G'_n = \prod_{j\neq n} G_j$. Then, for all $z\in Z_n$,
$$d(\gamma z, \gamma_n z) = d(\gamma'_n z, z) = d(\gamma'_n z_0, z_0)$$
since $z\in Z_n$ and we have a product action on $Z$. Thus, for each $\gamma\in\Gamma$, the distance between $\gamma\psi(z) = \psi(\gamma z)$ and $\psi(\gamma_n z)\in \psi(Z_n)$ is bounded independently of $z\in Z_n$. Since the action is reduced, it follows that $\psi(Z_n) = X$. At this point, in order to conclude the proof and exhibit $\wt\tau$ via $\psi|_{Z_n}:Z_n\cong X$, it only remains to show that $G_j$ fixes $Z_n$ (equivalently, that $Z_j$ is trivial) for all $j\neq n$. Indeed, the continuity of $\wt\tau$ as defined in Section~\ref{sec_group_action} follows from the continuity of $G$ on $Z$, and then $\wt\tau$ must range in $H$ by irreducibility of $\Gamma$.

To that end, we note that the above discussion did not depend on the choice of the ``copy'' $Z_n$; hence $\psi(gZ_n)=X$ also holds for all $g\in G$. Now pick $g_j\in G_j$ for some $j\neq n$. Since $\psi|_{Z_n}$ and $\psi|_{g_j Z_n}$ are isometries onto $X$, the map $T:X\to X$ defined by $\psi\circ g_j \circ(\psi|_{Z_n})^{-1}$ is an isometry. Since $(\psi|_{Z_n})^{-1}(x)$ and $g_j (\psi|_{Z_n})^{-1}(x)$ belong to a copy of $Z_j$, we have
$$d(T(x), x) = d(g_j(\psi|_{Z_n})^{-1}(x),(\psi|_{Z_n})^{-1}(x)) = d(g_j z_0, z_0)$$
independently of $x\in X$. Thus $T$ is a Clifford translation; since we assume that $X$ is not isometric to a Hilbert space, $T$ is trivial by~(iv) in Lemma~\ref{lemma_reduced}. This proves that $g_j$ fixes $z_0$ and therefore also $Z_n$, finishing the proof of Theorem~\ref{thm_non-proper} in this case.
\end{proof}

\subsection{}
\label{sec_Rd}%
We still have to prove Theorem~\ref{thm_non-proper} in the case where $X$ is a Hilbert space. We provide a proof just for the sake of completeness; this is a very special case for which our methods are rather coarse~--- indeed it is not so natural to assume the action reduced here since the linear structure allows for stronger results, see~\cite{Shalom00}. Moreover, Theorem~\ref{thm_non-proper} does not hold as stated when $X={\bf R}^d$ even in the simplest cases:

\begin{exo}
Consider the affine groups $G_1 = G_2 = {\bf R}\rtimes\{\pm 1\}$ and set $G=G_1\times G_2$. Consider the irreducible unifom lattice $\Gamma<G$ defined by
$$\Gamma\ \defq\ \Big\{\big((n+m\sqrt{2}, \epsi);(n-m\sqrt{2}, \epsi)\big) : n,m\in {\bf Z}, \epsi=\pm 1\Big\}.$$
The $\Gamma$\ti action on $X={\bf R}$ defined by $((n+m\sqrt{2}, \epsi);(n-m\sqrt{2}, \epsi)\big)x = \epsi x + n$ is unbounded and reduced. However it does not extend continuously to $G$. Indeed, if it did, the linear part of the $G$\ti action would have to factor through one of the $G_i$, say $G_1$, because the two canonical projections are the only homomorphisms $\{\pm1\}^2\to\mathbf{O}(1)=\{\pm1\}$ that are non-trivial on the diagonal. It then follows that $G_2$ acts trivially altogether, which is impossible since the homomorphism $n+m\sqrt{2}\mapsto n$ does not extend continuously to ${\bf R}$.

A yet simpler example is $G=({\bf Z}\rtimes\{\pm1\})^2$, $\Gamma={\bf Z}^2\rtimes \{\pm1\}$. Similar examples can be constructed in higher dimension.
\end{exo}

The above example is rather typical for counter-examples to the statement of Theorem~\ref{thm_non-proper}:

\begin{thm}
\label{thm_counter}%
The statement of Theorem~\ref{thm_non-proper} holds for any complete \cat0 space $X$ unless (i)~$X$ is isometric to ${\bf R}^d$ for some $d\in{\bf N}$ and (ii)~the linear part of the $\Gamma$\ti action is irreducible and is the restriction of a $G_i$\ti subrepresentation of $L^2(G/\Gamma)$ endowed with the quasi-regular $G_i$\ti representation. Moreover,~(ii) holds for at least two distinct indices~$i$.
\end{thm}

In particular, this theorem completes the proof of Theorem~\ref{thm_non-proper} since it is assumed in the latter that $X\neq{\bf R}^d$.

\smallskip

We begin with a general observation about lattices in products:

\begin{lemma}
\label{lem_two_ext}%
Let $\Gamma$ be a lattice in a product $G=G_1\times \cdots\times G_n$ of locally compact groups and let $\tau:\Gamma\to H$ be any homomorphism to a topological group $H$. Suppose that for two distinct indices $i$ there are continuous homomorphisms $\tau_i:G_i\to H$ such that the composed homomorphisms $\wt\tau_i: G\twoheadrightarrow G_i\to H$ extend $\tau$. Then $\tau(\Gamma)$ is relatively compact in $H$.
\end{lemma}

\begin{proof}
Let $i\neq i'$ be the two indices and define a continuous map $\sigma:G\to H$ by $\sigma(g) = \wt\tau_i(g)\wt\tau_{i'}(g)^{-1}$. Since $\sigma$ descends to a map $G/\Gamma\to H$ we obtain a probability Radon measure $\mu$ on $H$ from the (normalized) invariant measure on $G/\Gamma$. Since $\sigma(g_i g) = \tau_i(g_i)\sigma(g)$ for all $g_i\in G_i$ and $g\in G$, the measure $\mu$ is invariant under $\tau_i(G_i) = \wt\tau_i(G)$, hence under $\tau(\Gamma)$. It remains to observe that the stabilizer $K<H$ of any Radon probability measure on $H$ is compact.  Since $\mu$ is Radon there is a compact set $C\se H$ with $\mu(C)>1/2$. Therefore $kC\cap C\neq\varnothing$ for all $k\in K$. It follows $K\se C C^{-1}$.
\end{proof}

(The proof is slightly shorter when $\Gamma$ is cocompact.)

\begin{proof}[Proof of Theorem~\ref{thm_counter}]
In view of the proof of Theorem~\ref{thm_non-proper} given in Section~\ref{sec_non-proper} under the assumption that $X$ was not a Hilbert space, we may now assume that $X$ is a Hilbert space. Recall that a $\Gamma$\ti action by isometries on $X$ is given by an orthogonal representation $\pi\to\mathbf{O}(X)$ and a cocycle $b:\Gamma\to X$. We claim that $\pi$ is irreducible. Indeed, if $X=X'\oplus X''$ were a non-trivial orthogonal decomposition preserved by $\pi$, then $Y=0 \oplus X''$ would contradict the condition of Definition~\ref{defi_reduced} since $\tau(\gamma)Y$ is at finite distance of $\pi(\gamma)Y = Y$ for all $\gamma\in\Gamma$.

We adopt now the notations of the proof of Theorem~\ref{thm_non-proper} given in Section~\ref{sec_non-proper}; the non-Hilbertian assumption on $X$ was only used at the very end in order to prove that $G_j$ fixes $Z_n$ (equivalently, that $Z_j$ is trivial) for all $j\neq n$. Therefore, we assume now that $Z_j$ is non-trivial for some $j\neq n$ and need to characterize $X$ and its $\Gamma$\ti action as in Theorem~\ref{thm_counter}. We obtained in the proof of Section~\ref{sec_non-proper} an isometry $\psi|_{Z_n}:Z_n\cong X$ such that for each $\gamma\in\Gamma$, the distance between $\gamma\psi(z)$ and $\psi(\gamma_n z)$ is bounded independently of $z\in Z_n$. Therefore, the new $\Gamma$\ti action on $X$ transported via $\psi$ (and $\Gamma\to G_n$) from the $G_n$\ti action on $Z_n$ (i.e. $\gamma x = \psi(\gamma_n \psi|_{Z_n}^{-1}(x))$) differs from the original $\Gamma$\ti action by its translation cocycle only. In other words, $\pi$ extends continuously to a homomorphism
$$\wt\pi_n: G \twoheadrightarrow G_n \xrightarrow{\ \pi_n\ } \mathbf{O}(X)$$
where the orthogonal group $\mathbf{O}(X)$ is endowed with the strong operator topology. Applying the same argument to $j\neq n$ we are in position to use Lemma~\ref{lem_two_ext} and conclude that $\pi$ ranges in a compact subgroup of $\mathbf{O}(X)$. The Peter-Weyl theorem implies that $X$ is finite-dimensional since $\pi$ is irreducible. It remains only to prove that $\pi_i\se L^2(G/\Gamma)|_{G_i}$ for $i=j,n$. The $G$\ti representation $\ro=\mathrm{Ind}_\Gamma^G\pi$ induced from $\pi$ is canonically isomorphic to $\wt\pi_n\otimes L^2(G/\Gamma)$ since $\pi= \wt\pi_n|_\Gamma$ and hence $\pi_n \se \ro|_{G_n}$. Likewise, $\ro \cong \wt\pi_j\otimes L^2(G/\Gamma)$ and hence $\ro|_{G_n}\cong \mathrm{dim}(\pi)L^2(G/\Gamma)$. Since $\pi_n$ is irreducible, the claim follows for $i=n$ and is esteablished in the same way for $i=j$.
\end{proof}

\section*{Appendix~A: Commensurator Superrigidity}
%
\renewcommand{\thethmapp}{A\arabic{thmapp}} 
Let $G$ be a locally compact $\sigma$\ti compact group, $\Gamma<G$ a cocompact (or square-integrable, weakly cocompact) lattice and $\Gamma<\Lambda<G$ a dense subgroup \emph{commensurating} $\Gamma$, \emph{i.e.} $\Gamma\cap\lambda\Gamma\lambda^{-1}$ has finite index in $\Gamma$ for all $\lambda\in\Lambda$; equivalently, all $\Gamma$\ti orbits in $\Lambda/\Gamma$ are finite. In an unpublished manuscript~\cite{MargulisCOM} (see also~\cite{BurgerICM}), Margulis proves the following theorem (under a more relaxed non-positive curvature assumption, assuming $G$ compactly generated, $\Gamma$ cocompact, and assuming there are no parallel orbits).

\begin{thmapp}
Suppose $\Lambda$ acts by isometries on a complete \cat0 space $X$ such that the resulting $\Gamma$\ti action is non-evanescent. Then, upon possibly passing to a non-empty $\Gamma$\ti invariant closed convex subspace, the $\Gamma$\ti action extends continuously to a $G$\ti action.
\end{thmapp}

Margulis' proof uses generalized harmonic maps. We give an elementary proof illustrating the techniques introduced above; in spirit, this is a non-linear analogue of~\cite{Shalom00}. In the particular case where both $G$ and $\isom(X)$ are simple algebraic groups, this result leads to Margulis' arithmeticity criterion, see~\cite{Margulis},\cite{A'Campo-Burger}.

\begin{proof}
Using Lemma~\ref{lemma_KK}, we may assume $G$ second countable; we can assume $\Lambda$ countable and thus $X$ separable upon passing to the closed convex hull of a $\Lambda$\ti orbit. Let $Y$ be the induced $G$\ti space, non-evanescent by Theorem~\ref{thm_ind_eva} (resp. Appendix~B). It is enough to show that there is a non-empty $G$\ti invariant closed convex subspace $Z\se Y$ such that for all $f,f'\in Z$ the function $d(f, f')$ on $G$ is right $\Lambda$\ti invariant, since then it is essentially constant and we get a $\Gamma$\ti equivariant isometric map $Z\to X$ by evaluation(s) (compare~\ref{sec_super_intermediate}). Let $\A$ be the net of finite $\Gamma$\ti invariant sets $\varnothing\neq A\se\Lambda/\Gamma$. For any $f\in Y$ and a.e. $g\in G$ let $F_A f (g)$ be the unique $x\in X$ minimising $\sum_{a\in A} d^2(af(g a), x)$; this is the barycentre construction (Section~\ref{sec_bary}) for the uniform measure on $A$. We thus obtain a well-defined $G$\ti equivariant map $F_A:Y\to Y$. The barycentre inequality~(\ref{eq_bary}) of Section~\ref{sec_bary} yields
\begin{multline}
\label{eq_F_contr}%
d^2(F_A f(g), F_A f'(g))\ \leq\ |A|^{-1} \sum_{a\in A} d^2(f(ga), f'(ga))\\
 - |A|^{-1} \sum_{a\in A} \Big( d(f(ga), f'(ga)) - d(F_A f(g), F_A f'(g)) \Big)^2
\end{multline}
The first term already implies that $F_A$ is non-expanding (by integrating over $G/\bigcap_{a\in A} a\Gamma a^{-1}$); likewise,~(\ref{eq_F_contr}) shows that if $f,f'$ are $F_A$\ti fixed, then $d(f, f')$ is invariant under the group $\Lambda_A$ generated by the preimage of $A$ in $\Lambda$. Every $F_A$\ti orbit being evanescent, it is bounded; thus $Y^{F_A}\neq\varnothing$ by Remark~\ref{rem_non-expanding} --~\emph{and we are done if $\Lambda$ is finitely generated}, taking $A$ large enough and $Z=Y^{F_A}$.

For general $\Lambda$, let $T_A$ be the set of $G$\ti components $C\in T$ (Remarks~\ref{rems_exi_bis}) such that $d(f, f')$ is $\Lambda_A$\ti invariant $\forall\,f,f'\in C$. It follows from the preceding that $T_A\neq \varnothing$; $T_A$ is convex (use \emph{e.g.} Proposition~\ref{prop_paral_proj}) and closed since $L^2$\ti convergence of functions implies a.e. subconvergence. Since $T$ is bounded by non-evanescence, the directed family $\{T_A\}_{A\in \A}$ has non-empty intersection by Theorem~\ref{thm_compact} and yields a component $Z$ as sought.
\end{proof}

\section*{Appendix~B: Induction for Certain Non-Uniform Lattices}
%
\setcounter{thmapp}{0}
\renewcommand{\thethmapp}{B\arabic{thmapp}} 
This appendix discusses non-uniform lattices. The cocompactness assumption was \emph{only} used in Sections~\ref{sec_ind_G} and~\ref{sec_ind_eva}, and it was needed only in defining the action on the induced space (Lemma~\ref{lemma_ind_cont}) and in proving Theorem~\ref{thm_ind_eva}. Thus, replacing them respectively with Lemma~\ref{lemma_ind_app} and Theorem~\ref{thm_ind_app} below, we conclude that all our results (including Appendix~A) hold as claimed in Theorem~\ref{thm_non_uniform}. We insist however that for classical lattices the integrability condition discussed below is dependent on Margulis' arithmeticity theorem (interestingly, Margulis originally proved arithmeticity of non-uniform lattices without~-- and before~-- his superrigidity).

Let $G$ be a locally compact second countable group, $\Gamma<G$ a lattice. We use the notation of Sections~\ref{sec_ind_G} and~\ref{sec_ind_eva}, \emph{e.g.} for $\chi:G\to\Gamma$, $\psi_x(g) \defq \chi(g)x$ and normalizing the covolume of $\Gamma$ to one.

\begin{defiapp}[{See~\cite[1.II]{Shalom00}}]
\label{def_sq_sum}%
The lattice $\Gamma$ is \emph{square-integrable} if it is finitely generated and if, for the length function $\ell$ associated to some (or equivalently any) finite generating set, there is a fundamental domain $\F\se G$ (with null boundary) such that
$$\int_\F \ell(\chi(g^{-1} h))^2\d h \ <\infty\kern1cm\forall\,g\in G.$$
(We note that when dealing with uniform lattices we never imposed finite generation, thus allowing for lattices in groups that are not compactly generated.)
\end{defiapp}

Y.~Shalom explains in~\cite[\S2]{Shalom00} why the condition of Definition~\ref{def_sq_sum} always holds for lattices as in Theorem~\ref{thm_Margulis}; B.~R\'emy proves in~\cite{Remy04} that it holds for all Kac-Moody lattices. We refer to~\cite{RemyAST},\cite{RemySURV} for general Kac-Moody groups, in particular for the following result of R\'emy: Any Kac-Moody group over ${\bf F}_q$ is an irreducible lattice in the product of its associated twin building groups (modulo its finite centre), when $q$ is large enough.

\smallskip

The following parallels~\cite[1.II]{Shalom00}:

\begin{lemmaapp}
\label{lemma_ind_app}%
Let $\F$ be as in Definition~\ref{def_sq_sum}. Then Lemma~\ref{lemma_ind_cont} and formula~(\ref{eq_act}) provide a well-defined continuous $G$\ti action by isometries on $L^{[2]}(G,X)^\Gamma \cong L^2(\F,X)$.
\end{lemmaapp}

\begin{proof}
The only additional verification we need to do is that for $f\in L^2(\F,X)$, $x\in X$ and $g\in G$ the integral $\int_\F d^2((gf)(h),x)\d h$ is finite. The latter is $\int_\F d^2(f(g^{-1}.h), \chi(g^{-1} h)^{-1}x)\d h$ in view of~(\ref{eq_act}); since $\int_\F d^2(f(g^{-1}.h),x)\d h = \int_\F d^2(f(h),x)\d h$ is finite, it is enough to show that $\int_\F d^2(\chi(g^{-1}h)^{-1} x, x)\d h$ is finite. Let $S\se \Gamma$ be a finite generating set and $\ell$ the associated length function; since $d(\chi(g^{-1} h)^{-1} x, x)$ is bounded by $\ell(\chi(g^{-1} h))\sup_{s\in S}d(s x, x)$, we conclude by square-integrability of $\Gamma$.
\end{proof}

The evanescence question (Theorem~\ref{thm_ind_eva}) is more difficult; we shall establish a geometric generalization of an argument given in the linear setting by Margulis~\cite[III.1]{Margulis}.

\begin{defiapp}[{\cite[III.1.8]{Margulis}}]
The lattice $\Gamma$ is \emph{weakly cocompact} if the $G$\ti representation $L^2_0(G/\Gamma)$ (\emph{i.e.} the orthogonal complement of the trivial representation in $L^2(G/\Gamma)$) does not almost have non-zero invariant vectors (compare~\ref{sec_Kazhdan}).
\end{defiapp}

The definition is obviously satisfied whenever $G$, or equivalently $\Gamma$, has Kazhdan's property~(T); this disposes right away with most higher rank groups. According to Margulis~\cite[II.1.12]{Margulis}, it also holds for connected semisimple Lie groups $G$ even when they are not Kazhdan, see~\cite{Bekka98} for a proof. Any Kac-Moody group over ${\bf F}_q$ whose Cartan matrix has finite entries is Kazhdan whenever~$q$ is large enough by a general result of Dymara-Januszkiewicz~\cite{Dymara-J02}.

\begin{remapp}
\label{rem_asympt_app}%
If there is an asymptotically invariant sequence $\{v_n\}$ (see Remark~\ref{rem_asympt_inv}) of non-negative functions in $L^2(G/\Gamma)$ such that for every relatively compact $C\se \F$ the integral $\int_C v_n$ tends to zero, then $\Gamma$ is not weakly cocompact. Indeed, in that case $\int_\F v_n\to 0$, thus the norm of the projection of $v_n$ on $L^2_0(G/\Gamma)$ tends to one, yielding (after renormalization) an asymptotically invariant sequence in $L^2_0$.
\end{remapp}

\begin{thmapp}
\label{thm_ind_app}%
Assume that $\Gamma$ is square-integrable and weakly cocompact. Let $X$ be a complete separable \cat0 space with a non-evanescent $\Gamma$\ti action by isometries. Then the $G$\ti action on $L^{[2]}(G,X)^\Gamma$ is non-evanescent.
\end{thmapp}

\begin{proof}
Fix $x_0\in X$. Suppose for a contradiction that there is an evanescent sequence $\{f_n\}$ in $Y=L^{[2]}(G,X)^\Gamma$ such that $d(f_n, \psi_{x_0})\to \infty$. Let $\eta$ be a non-negative continuous function on $G$ of integral one; we may assume that $\eta$ has compact support $K$ with $K\setminus \F$ null. For $g\in G$, we want to define $\bar f_n(g)\in X$ as the barycentre of $h\mapsto f_n(h^{-1}g)$ with respect to the measure $\eta(h)\d h$; we thus have to prove that for (some, hence any $x\in X$) the integral $\int_\F d^2(f_n(h^{-1}g), x)\eta(h)\d h$ is finite for a.e. $g\in \F$ (hence a.e. $g\in G$). This follows from Tonelli's theorem applied to
\begin{multline}
\label{eq_Tonelli}%
\int_\F\int_\F d^2(f_n(h^{-1}g), x)\eta(h)\d h \d g\ =\ \int_\F \eta(h) \int_\F d^2(f_n(h^{-1}g), x)\d g \d h\\
=\ \int_K \eta(h) d^2(hf_n, \psi_x)\d h\ \leq\ \sup \big\{ d^2(hf_n, \psi_x) : h\in K \big\}\ < \infty.
\end{multline}
By the definition of barycentres (Section~\ref{sec_bary}),
\beq
\label{eq_bar_f_f}%
d^2(\bar f_n(g), x) \ \leq\ \int_\F d^2(f_n(h^{-1}g), x)\eta(h)\d h\kern1cm\forall\,x\in X.
\eeq
Thus, $\bar f_n$ is square-integrable because we apply~(\ref{eq_Tonelli}) to
$$\int_\F d^2(\bar f_n(g), x)\d g \ \leq\ \int_\F\int_\F d^2(f_n(h^{-1}g), x)\eta(h)\d h \d g.$$
Since in addition $\bar f_n$ is $\Gamma$\ti equivariant by definition, $\bar f_n\in Y$. Setting $x=f_n(g)$ in~(\ref{eq_bar_f_f}) yields
\beq
\label{eq_f_fbar}%
d^2(\bar f_n, f_n) \ \leq\ \int_\F\int_\F d^2(f_n(h^{-1}g), f_n(g))\eta(h)\d h \d g \ =\ \int_K \eta(h) d^2(h f_n, f_n)\d h
\eeq
which is bounded independently of $n$ by evanescence of $\{f_n\}$. It follows that $\{\bar f_n\}$ is also an evanescent sequence with $d(\bar f_n, \psi_{x_0})\to \infty$. Define now $\fhi_n\in L^2(G/\Gamma)$ by $\fhi_n(h) \defq d(\bar f_n(h), \chi(h)x_0)$. We claim that $\{\fhi_n\}$ is an evanescent sequence in the linear $G$\ti space $L^2(G/\Gamma)$. Indeed, since $(g\fhi_n)(h) = d((g\bar f_n)(h), \chi(g^{-1} h) x_0)$, the triangle inequality gives
$$\big| (g\fhi_n)(h) - \fhi_n(h) \big| \ \leq\ d((g\bar f_n(h), \bar f_n(h)) + d(\chi(g^{-1} h) x_0, x_0)$$
so that by Minkowski's inequality and the definition of $\psi_{x_0}$
$$\|g\fhi_n - \fhi_n\| \ \leq\ d(g \bar f_n, \bar f_n) + d(g\psi_{x_0}, \psi_{x_0}).$$
The first term is bounded over compact sets by evanescence of $\{\bar f_n\}$ and the second by continuity of the $G$\ti action on $Y$; the claim follows. In particular, since $\|\fhi_n\| = d(\bar f_n, \psi_{x_0})\to\infty$, the sequence $\{v_n = \fhi_n / \|\fhi_n\|\}$ is asymptotically invariant.

The goal now is to contradict weak cocompactness by applying Remark~\ref{rem_asympt_app} to $\{v_n\}$. Therefore it suffices to show that for any relatively compact $C\se \F$ the integral $\int_C \fhi_n(g) \d g$ is bounded independently of $n$.

Since the $\Gamma$\ti action on $X$ is not weakly evanescent (Proposition~\ref{prop_ceva}), there is by Lemma~\ref{lemma_eva_lin} a finite set $F\se \Gamma$, $\lambda>0$ and $d_0\geq 0$ such that $\sup_{\gamma\in F} d(\gamma^{-1} x, x) \geq \lambda d(x,x_0) - d_0$ for all $x\in X$. Thus,
$$\lambda \fhi_n(g) \ \leq\ \sup_{\gamma\in F} d(\gamma^{-1} \bar f_n(g), \bar f_n(g)) + d_0 \ \leq\ \sum_{\gamma\in F} d(\bar f_n(g\gamma), \bar f_n(g)) + d_0 \kern1cm(\mathrm{a.e.}\ g\in\F).$$
Therefore,
\begin{multline}
\label{eq_integral_eva}%
\lambda\int_C \fhi_n(g) \d g \ \leq\ \sum_{\gamma\in F} \int_C d(\bar f_n(g\gamma), \bar f_n(g))\d g + d_0\\
\leq\ \sum_{\gamma\in F} \left(\int_C d^2(\bar f_n(g\gamma), \bar f_n(g))\d g\right)^{1/2} + d_0.
\end{multline}
It is now enough to prove that for all $\gamma\in F$ the integral $\int_C d^2(\bar f_n(g\gamma), \bar f_n(g))\d g$ is bounded independently of $n$. To that end, set $K_{\gamma,g} \defq g\gamma^{-1} g^{-1} K$, $\eta_{\gamma, g}(\cdot) \defq \eta(g\gamma g^{-1}\cdot)$ for any $\gamma\in F$, $g\in G$. By definition of $\bar f_n$ and change of variable, $\bar f_n(g\gamma)$ is the minimiser $y\in X$ of $\int_{K_{\gamma,g}} d^2( f_n(h^{-1} g), y) \eta_{\gamma, g}(h) \d h$. Therefore, the inequality corresponding to~(\ref{eq_bar_f_f}) yields
$$d^2(\bar f_n(g\gamma), \bar f_n(g)) \ \leq\ \int_{K_{\gamma,g}} d^2( f_n(h^{-1} g), \bar f_n(g)) \eta_{\gamma, g}(h) \d h \ \leq\ \|\eta\|_\infty \int_{K_{\gamma,g}} d^2( f_n(h^{-1} g), \bar f_n(g)) \d h.$$
Thus,
$$\int_C d^2(\bar f_n(g\gamma), \bar f_n(g)) \d g \ \leq\ \|\eta\|_\infty \int_C \int_L d^2( f_n(h^{-1} g), \bar f_n(g)) \d h \d g,$$
where the integral over $L \defq \bigcup\{ K_{\gamma, g} : g\in C\}$ is finite because we bound the above double integral by
$$\int_L \int_C d^2((h f_n)(g), \bar f_n(g)) \d g \d h \ \leq\ \int_L d^2(h f_n, \bar f_n) \d h$$
which is finite by relative compactness of $L$. Moreover, the latter term is bounded independently of $n$ in view of
$$d(h f_n, \bar f_n) \leq d(h f_n, f_n) + d(f_n, \bar f_n)$$
since the first summand here is bounded by evanescence and the second has been treated previously with~(\ref{eq_f_fbar}). This concludes the proof.
\end{proof}

%
\ifx\undefined\bysame
\newcommand{\bysame}{\leavevmode\hbox to3em{\hrulefill}\,}
\fi

\end{document}